\documentclass[11pt,a4paper]{amsart}
\usepackage[T1]{fontenc}
\usepackage[latin1]{inputenc}
\usepackage{amssymb}

\makeatletter
  \theoremstyle{plain}
  \newtheorem{thm}{Theorem}
\newtheorem{appthm}{Theorem}
  \theoremstyle{remark}
  \newtheorem*{rem*}{Remark}
  \theoremstyle{remark}
  \newtheorem*{acknowledgement*}{Acknowledgement}
 \theoremstyle{definition}
 \newtheorem*{defn*}{Definition}
  \theoremstyle{plain}
  \newtheorem{prop}{Proposition}
  \theoremstyle{remark}
  \newtheorem{rem}{Remark}
  \theoremstyle{plain}
  \newtheorem{lem}{Lemma}
\newtheorem{applem}{Lemma}
  \theoremstyle{remark}
  \newtheorem*{note*}{Note}

\usepackage[english]{babel}
\makeatother
\begin{document}
\newcommand{\fs}{f_{1},\ldots f_{r}}
\newcommand{\Sing}{\mathrm{Sing}}
\newcommand{\Spec}{\mathrm{Spec}\ }
\newcommand{\sB}{\mathsf{B}}
\newcommand{\Proj}{\mathrm{Proj}\ }

\title{The density of integral points on complete intersections}

\author{Oscar Marmon}
\address{Department of Mathematical Sciences\\ Chalmers University of Technology and Göteborg University\\ SE-412 96 Göteborg\\ Sweden}
\email{marmon@math.chalmers.se, salberg@math.chalmers.se}

\maketitle
\let\languagename\relax

\vspace{-0.5cm}
\begin{center}
\small\textsc{with an appendix by Per Salberger}
\end{center}

\begin{abstract}
In this paper, an upper bound for the number of integral points of
bounded height on an affine complete intersection defined over $\mathbb{Z}$
is proven. The proof uses an extension to complete intersections of
the method used for hypersurfaces by Heath-Brown \cite{Heath-Brown},
the so called {}``$q$-analogue'' of van der Corput's AB process.
\end{abstract}

\section{Introduction}

If $X$ is an affine algebraic set defined by a set of equations \[
f_{i}(x_{1},\ldots x_{n})=0,i=1,\ldots,r\]
with integral coefficients, and if $\sB$ is a box in $\mathbb{R}^{n}$
- that is, a product of closed intervals - then we define the quantity\[
N(X,\sB)=\#\left\{ \mathbf{x}=(x_{1},\ldots,x_{n})\in\mathbb{Z}^{n};f_{i}(\mathbf{x})=0,\mathbf{x}\in\sB\right\} .\]
If $m$ is a positive integer, and if \textbf{$\sB$} is small enough
as to contain at most one representative of each congruence class
modulo $m$, then we define\[
N(X,\sB,m)=\#\left\{ \mathbf{x}=(x_{1},\ldots,x_{n})\in\mathbb{Z}^{n};f_{i}(\mathbf{x})\equiv0\pmod m,\mathbf{x}\in\sB\right\} .\]
Since $N(X,\sB)\leq N(X,\sB,m)$ one can obtain upper bounds for $N(X,\sB)$
by considering $N(X,\sB,m)$ for suitably chosen $m$. If $\sB=[-B,B]^{n}$
for some $B>0$ we write \[
N(X,B)=N(X,\mathsf{B})\text{ and }N(X,B,m)=N(X,\mathsf{B},m).\]
Throughout this paper we shall be concerned with the case when $X$
is a complete intersection, that is, when $\dim X=n-r$, where $r$
is the number of equations defining $X$ in $\mathbb{A}^{n}$. Our
main concern shall be to find an upper bound for $N(X,B)$. One result
in this direction is the following, by Fujiwara \cite{Fujiwara88}:
let $X$ be a non-singular hypersurface in $\mathbb{A}^{n}$ defined
by the vanishing of a polynomial $f$ with integer coefficients, of
degree at least $2$. Then $N(X,B)\ll_{f,n}B^{n-2+2/n}$ for $n\geq4$.
Fujiwara proved this by exhibiting an asymptotic formula for $N(X,B,p)$
for primes $p$, the proof of which uses the estimates for exponential
sums by Deligne \cite{Deligne} as a key tool. Heath-Brown \cite{Heath-Brown}
was able to sharpen the exponent to $n-2+2/(n+1)$ by averaging over
primes in an interval. In the same paper he introduced a new technique,
the so called $q$-analogue of van der Corput's method. He could then
prove the bound\begin{equation}
N(X,B)\ll_{f,n}B^{n-3+15/(n+5)}\label{eq:H-B}\end{equation}
for a non-singular hypersurface $X$ defined by a polynomial $f$
of degree at least 3 (Theorem 2 in \cite{Heath-Brown}), by considering
$N(X,B,pq)$ for two suitable primes $p$ and $q$. 

In this paper we will generalize the method of Heath-Brown to complete
intersections of arbitrary codimension. We shall use the following
notation: if $X$ is a scheme over $\mathbb{Z}$ we let $X_{\mathbb{Q}}=X\times_{\Spec\mathbb{Z}}\mathbb{Q}$
and $X_{q}=X_{\mathbb{F}_{q}}=X\times_{\Spec\mathbb{Z}}\mathbb{F}_{q}$
for every prime $q$.

\begin{thm}
\label{thm:uniform}Let\[
X=\Spec\mathbb{Z}[X_{1},\ldots,X_{n}]/(f_{1},\ldots,f_{r}),\]
where the leading forms $F_{1},\ldots,F_{r}$ of $f_{1},\ldots,f_{r}$
are of degree $\geq3$, and let \[
Z=\Proj\mathbb{Z}[X_{1},\ldots,X_{n}]/(F_{1},\ldots,F_{r}).\]
Assume that $Z_{\mathbb{Q}}$ is non-singular of codimension $r$
in $\mathbb{P}_{\mathbb{Q}}^{n-1}$. Then, if $n\geq4r+2$, we have
for $B\geq1$\[
N(X,B)\ll_{n,d,\epsilon}B^{n-3r+r^{2}\frac{13n-5-3r}{n^{2}+4nr-n-r-r^{2}}}(\log B)^{n/2}\left(\sum_{i=1}^{r}\log\left\Vert F_{i}\right\Vert \right)^{2r+1},\]
where $d=\max_{i}(\deg f_{i}).$
\end{thm}
\begin{rem*}
The factor $(\log B)^{n/2}$ can in fact be disposed of, and we sketch
in the end of Section \ref{sec:Proof-of-the-Main-Result} how this
can be done.
\end{rem*}
The estimate given by Theorem $\ref{thm:uniform}$ in the case $r=1$
is in fact slightly sharper than (\ref{eq:H-B}), owing to the use
of estimates by Katz \cite{Katz} on exponential sums modulo $q$.
Theorem \ref{thm:uniform} is a corollary to the following theorem. 

\begin{thm}
\label{thm:main}Let\[
X=\Spec\mathbb{Z}[X_{1},\ldots,X_{n}]/(f_{1},\ldots,f_{r}),\]
where $r<n$ and the leading forms $F_{1},\ldots,F_{r}$ of $f_{1},\ldots,f_{r}$
are of degree $\geq3$, and let \[
Z=\Proj\mathbb{Z}[X_{1},\ldots,X_{n}]/(F_{1},\ldots,F_{r}).\]
Let $B$ be a positive number, and let $p$ and $q$ be primes, with
$2p<2B+1<q-p$, such that both $Z_{p}$ and $Z_{q}$ are non-singular
of dimension $n-1-r$. Then we have\begin{multline*}
N(X,B,pq)=\frac{(2B+1)^{n}}{p^{r}q^{r}}+O_{n,d}\left(B^{(n+1)/2}p^{-r/2}q^{(n-r-1)/4}(\log q)^{n/2}\right.\\
+B^{(n+1)/2}p^{(n-2r)/2}q^{-1/4}(\log q)^{n/2}+B^{n/2}p^{-r/2}q^{(n-r)/4}(\log q)^{n/2}\\
\left.+B^{n/2}p^{(n-r)/2}(\log q)^{n/2}+B^{n}p^{-(n+r-1)/2}q^{-r}+B^{n-1}p^{-r+1}q^{-r}\right),\end{multline*}
where $d=\max_{i}(\deg f_{i}).$
\end{thm}
The proof of Theorem \ref{thm:main} is carried out in Section \ref{sec:Proof-of-the-Main-Result}
and more or less follows \cite{Heath-Brown}. However, in contrast
to Heath-Brown, we do not use Poisson summation, but a more direct
approach.

We also prove, in Section \ref{sec:modq}, a generalization (and slight
sharpening) of Theorem 3 in \cite{Heath-Brown}, a weighted asymptotic
formula for the density of $\mathbb{F}_{q}$-points on affine complete
intersections defined over $\mathbb{F}_{q}$. However, for the proof
of Theorem \ref{thm:main}, we will use an unweighted version of this
result, proven by Salberger in an Appendix to this paper. This is
because we desire an unweighted asymptotic formula in Theorem \ref{thm:main}.

\begin{acknowledgement*}
I wish to thank my supervisor Per Salberger for introducing me to
the topic of this paper, and for numerous helpful suggestions during
the way. 
\end{acknowledgement*}

\section{Preliminary Results from Algebraic Geometry}

We recall some facts from algebraic geometry that will provide helpful
tools for proving our main results.

\begin{defn*}
Let $X$ be a scheme. A point $x\in X$ is a \emph{singular point
of $X$} if the local ring $\mathcal{O}_{X,x}$ is not a regular local
ring. $X$ is said to be \emph{singular} if it has singular points,
and \emph{non-singular} if not. We denote the \emph{singular locus
of $X$ -} the set of singular points \emph{}- by $\Sing X$.

If $X$ is a scheme and $x$ a point on $X$, then $\mathcal{O}_{x}$
is the local ring at $x$, $\mathfrak{m}_{x}$ its maximal ideal and
$\kappa(x)=\mathcal{O}_{x}/\mathfrak{m}_{x}$ the residue field of
$x$. If $X\to Y$ is a morphism of schemes, $\Omega_{X/Y}$ denotes
the sheaf of relative differentials of $X$ over $Y$, and we abbreviate
$\Omega_{X/\Spec R}=\Omega_{X/R}$. 
\end{defn*}
We have the following characterization of singular points on a scheme.

\begin{prop}
\label{pro:singular}Let $X$ be a scheme of finite type over a perfect
field $k.$ Suppose that $X$ is equidimensional of dimension $n$.
Then for every point $x\in X$, the following conditions are equivalent:
\end{prop}
\begin{itemize}
\item [(i)]$x$ is a singular point of $X$;
\item [(ii)]$\dim_{\kappa(x)}\Omega_{X/k,x}\otimes_{\mathcal{O}_{x}}\kappa(x)>n.$
\end{itemize}
\begin{proof}
Since this is a local question, we can assume that $X=\Spec R$ with
$R$ equidimensional. Suppose $x=\mathfrak{p}\in\Spec R$. Then we
have, by \cite[Ex. 14.36]{Sharp},\begin{equation}
\begin{aligned}n= & \mathrm{ht}\mathfrak{p}+\dim R/\mathfrak{p}\\
= & \dim\mathcal{O}_{x}+\mathrm{tr.d.}\kappa(x)/k.\end{aligned}
\label{eq:Sharp}\end{equation}
By definition, $x$ is singular if and only if\[
\dim_{\kappa(x)}\mathfrak{m}_{x}/\mathfrak{m}_{x}^{2}>\dim\mathcal{O}_{x}.\]
Furthermore, by \cite[Ex. II.8.1]{Hartshorne}, we have an exact sequence
of $\kappa(x)$-vector spaces \[
0\to\mathfrak{m}_{x}/\mathfrak{m}_{x}^{2}\to\Omega_{\mathcal{O}_{x}/k}\otimes_{\mathcal{O}_{x}}\kappa(x)\to\Omega_{\kappa(x)/k}\to0.\]
Since $\Omega_{\mathcal{O}_{x}/k}$ is equal to the stalk $\Omega_{X/k,x}$
of the sheaf of relative differentials, and since $\dim_{\kappa(x)}\Omega_{\kappa(x)/k}=\mathrm{tr.d.}\kappa(x)/k$
by \cite[Thm. II.8.6A]{Hartshorne}, this implies that\[
\dim_{\kappa(x)}\Omega_{X/k,x}\otimes_{\mathcal{O}_{x}}\kappa(x)=\dim_{\kappa(x)}\mathfrak{m}_{x}/\mathfrak{m}_{x}^{2}+\mathrm{tr.d.}\kappa(x)/k.\]
In view of (\ref{eq:Sharp}) it follows that $x\in\Sing X$ if and
only if\[
\dim_{\kappa(x)}\Omega_{X/k,x}\otimes_{\mathcal{O}_{x}}\kappa(x)>\dim\mathcal{O}_{x}+\mathrm{tr.d.}\kappa(x)/k=n.\]

\end{proof}
\begin{rem}
\label{rem:varphi(x)}By \cite[Ex. II.5.8]{Hartshorne} the function
\[
\varphi(x)=\dim_{\kappa(x)}\Omega_{X/k,x}\otimes_{\mathcal{O}_{x}}\kappa(x)\]
is upper semicontinuous, so that in the situation described in the
proposition, $\Sing X$ is a closed subscheme of $X$. 
\end{rem}
{}

\begin{rem}
\label{rem:smooth}The proposition also shows that for $X$ equidimensional
and of finite type over a perfect field $k$, $X$ is non-singular
if and only if it is \emph{smooth over $k$} (see \cite[Ch. III.10]{Hartshorne}).
\end{rem}
{}

\begin{rem}
The particular case where we will use the proposition is for $X$
a complete intersection of positive dimension in projective space
over a perfect field. Such $X$ are indeed equidimensional, since
firstly, any local complete intersection is Cohen-Macaulay (\cite[Prop. 8.23]{Hartshorne})
and thus locally equidimensional, and secondly, a complete intersection
in $\mathbb{P}_{k}^{n}$ of dimension $\geq1$ is connected (\cite[Ex. III.5.5]{Hartshorne}).
\end{rem}
When working in a projective space $\mathbb{P}^{n}$ with homogeneous
coordinates $x_{0},\ldots,x_{n}$ we denote by $\check{\mathbb{P}^{n}}$
the dual projective space with homogeneous coordinates $\xi_{0},\ldots,\xi_{n}.$
For a point $\mathbf{a}=(a_{0},\ldots,a_{n})$ in $\check{\mathbb{P}^{n}}$
we will let $H_{\mathbf{a}}$ denote the hyperplane defined in $\mathbb{P}^{n}$
by the equation $\mathbf{a}\cdot\mathbf{x}=a_{0}x_{0}+\ldots+a_{n}x_{n}=0$.
We begin by proving the following corollary to Bertini's Theorem.
By convention, the dimension of the empty set is defined to be $-1$.

\begin{lem}
\label{lem:Bertini.}Let $k$ be an algebraically closed field. Let
$X$ be a non-empty complete intersection in $\mathbb{P}_{k}^{n}$.
Suppose that\[
\dim\mathrm{Sing}X=s.\]
Then there is a hyperplane $H$ such that $\dim(X\cap H)=\dim X-1$
and \[
\dim\mathrm{Sing}(X\cap H)<\max(s,0).\]

\end{lem}
\begin{proof}
The case $s=-1$ follows immediately from Bertini's Theorem \cite[Cor 6.11(2)]{Jouanolou}.
($X$ is then smooth over $k$ by Remark \ref{rem:smooth}.) If $s\geq0$,
let $Y=X\setminus\Sing X$, so that $Y$ is smooth. Then, by Bertini's
Theorem, there exists a non-empty Zariski open subset $U$ of $\check{\mathbb{P}_{k}^{n}}$
such that for hyperplanes $H_{\mathbf{a}}$ parametrized by closed
$k$-points $\mathbf{a}$ in $U$, $Y\cap H_{\mathbf{a}}$ is smooth
and thus non-singular by Remark \ref{rem:smooth}. Hence, for $\mathbf{a}\in U(k)$
we have\begin{equation}
\Sing(X\cap H_{\mathbf{a}})\subseteq\Sing X\cap H_{\mathbf{a}}.\label{eq:Sing(Xcapa)}\end{equation}
Furthermore, there are non-empty open sets $U',U''$ such that for
all closed $k$-points $\mathbf{a}$ of $U'$, no irreducible component
of $\Sing X$ of dimension $s$ is contained in $H_{\mathbf{a}}$,
and for $\mathbf{a}\in U''(k)$ no irreducible component of $X$ is
contained in $H_{\mathbf{a}}$. Then we have, for $\mathbf{a}\in U\cap U'\cap U''(k)$,
that $\dim(X\cap H_{\mathbf{a}})=\dim X-1$ and $\dim\Sing(X\cap H_{\mathbf{a}})<s$.
\end{proof}
\begin{rem}
\label{rem:katzlemma3}For any hyperplane $H$ such that $\dim X\cap H=\dim X-1$,
$\dim\Sing(X\cap H)\geq\dim\mathrm{Sing}X-1$ (see \cite[Lemma 3]{Katz}).
\end{rem}
The next lemma is an {}``effective'' version of Bertini's Theorem.
For a more explicit result of the same type, see \cite{Ballico}.

\begin{lem}
\label{lem:Phi}Let $n,r,d_{1},\ldots,d_{r}$ be natural numbers,
and let $F_{1},\ldots,F_{r}$ be forms in $X_{0},\ldots,X_{n}$ with
integer coefficients, and with $\deg F_{i}=d_{i}$. Let $V=\Proj\mathbb{Z}[X_{0},\ldots,X_{n}]/(F_{1},\ldots,F_{r})$,
and suppose that $V_{\mathbb{Q}}$ has dimension $n-r\geq0$. Then
for every prime $q$ such that $V_{q}$ has dimension $n-r$, there
is a non-zero form $\Phi_{q}\in\mathbb{F}_{q}[\xi_{0},\ldots,\xi_{n}]$
with degree bounded in terms of $n$ and $d_{1},\ldots,d_{r}$ only,
such that for every point $\mathbf{a}=(a_{0},\ldots,a_{n})\in\check{\mathbb{P}_{\mathbb{F}_{q}}^{n}}$
satisfying $\Phi_{q}(a_{0},\ldots,a_{n})\neq0$ we have
\begin{itemize}
\item [(i)] $\dim\mathrm{Sing}(V_{q}\cap H_{\mathbf{a}})=\max(-1,\dim\mathrm{Sing}V_{q}-1)$
\item [(ii)]$\dim V_{q}\cap H_{\mathbf{a}}=\dim V_{q}-1$.
\end{itemize}
In particular, for each $q\geq q_{0}=q_{0}(n,d_{1},\ldots,d_{r})$
there is an $\mathbf{a}\in\check{\mathbb{P}_{\mathbb{F}_{q}}^{n}}$
with the properties (i) and (ii).
\end{lem}
\begin{proof}
We let $\mathbb{P}_{i}$, for each $i=1,\ldots,r$, be the projective
space over $\mathbb{Z}$ parametrizing all hypersurfaces in $\mathbb{P}_{\mathbb{Z}}^{n}$
of degree $d_{i}$ (as a Hilbert scheme), and work in the large multiprojective
space $\mathbf{P}=\mathbb{P}_{1}\times\ldots\times\mathbb{P}_{r}.$
For a $k$-point in $\mathbf{P}$ representing a tuple $(F_{1},\ldots,F_{r})$
we write $V(F_{1},\ldots,F_{r})$ for the intersection of the corresponding
$r$ hypersurfaces in $\mathbb{P}_{k}^{n}$. Let $W\subseteq\mathbf{P}\times\check{\mathbb{P}_{\mathbb{Z}}^{n}}\times\mathbb{P}_{\mathbb{Z}}^{n}$
be defined as the closed set of points $P\in\mathbf{P}\times\check{\mathbb{P}_{\mathbb{Z}}^{n}}\times\mathbb{P}_{\mathbb{Z}}^{n}$
representing $(F_{1},\ldots,F_{r},\mathbf{a},\mathbf{x})$ that satisfy\[
\mathbf{x}\in V(F_{1},\ldots,F_{r})\cap H_{\mathbf{a}}.\]
 Let\[
\pi:W\to\mathbf{P}':=\mathbf{P}\times\check{\mathbb{P}_{\mathbb{Z}}^{n}}\]
be the projection. The function $\varphi(P):=\dim_{\kappa(P)}\Omega_{W/\mathbf{P}',P}$
is upper semicontinuous (see Remark \ref{rem:varphi(x)}), so the
set \[
S=\left\{ P\in W;\varphi(P)\geq n-r\right\} \]
is closed. Now, let $\tilde{\pi}:S\to\mathbf{P}'$ be the restriction
of $\pi$ to $S$, and let for every $s\in\{-1,0,1,\ldots,n\}$\[
A_{s}=\left\{ Q\in\mathbf{P}';\dim\tilde{\pi}^{-1}(Q)\geq s\right\} .\]
By Chevalley's Semicontinuity Theorem \cite[Cor 13.1.5]{EGAIV(III)},
$A_{s}$ is closed in $\mathbf{P}'$, as is the set \[
D=\left\{ Q\in\mathbf{P}';\dim\pi^{-1}(Q)\geq n-r\right\} .\]
For each $s\in\{-1,0,\ldots,n\},$ let $T_{s}=D\cup A_{s}$. Then
$T_{s}$ is closed as well, so there exist multihomogeneous forms
$H_{1}^{s},\ldots,H_{t}^{s}$ over $\mathbb{Z}$ that define $T_{s}.$ 

For a closed $k$-point $P\in W$ representing $(F_{1},\ldots,F_{r},\mathbf{a},\mathbf{x})$
we have an isomorphism of stalks $\Omega_{W/\mathbf{P}',P}\cong\Omega_{Y/k,\mathbf{x}}$,
where \[
Y=V(F_{1},\ldots,F_{r})\cap H_{\mathbf{a}}\subseteq\mathbb{P}_{k}^{n}.\]
Thus, for each tuple $(F_{1},\ldots,F_{r},\mathbf{a})$ such that
both $V=V(F_{1},\ldots,F_{r})$ and $V\cap H_{\mathbf{a}}$ are complete
intersections of codimension $r$ and $r+1$, respectively, the fiber
$\tilde{\pi}^{-1}(F_{1},\ldots,F_{r},\mathbf{a})$ is precisely $\mathrm{Sing}(V\cap H_{\mathbf{a}})$
by Proposition \ref{pro:singular}. For every other point $(F_{1},\ldots,F_{r},\mathbf{a})$
we have $\tilde{\pi}^{-1}(F_{1},\ldots,F_{r},\mathbf{a})=\mathbb{P}_{k}^{n}$.
We conclude that $T_{s}$, for each $s$, is the set of tuples $(F_{1},\ldots,F_{r},\mathbf{a})$
such that $V(F_{1},\ldots,F_{r})\cap H_{\mathbf{a}}$ either has codimension
$\leq r$ or has a singular locus of dimension at least $s$. In particular,
if  we have a closed $k$-point $Q\in\mathbf{P}$ representing $(F_{1},\ldots,F_{r})$
such that $V=V(F_{1},\ldots,F_{r})$ satisfies \begin{equation}
\dim V=n-r,\quad\dim\mathrm{Sing}V=s,\label{eq:dimSing}\end{equation}
and if $\pi_{s}:T_{s}\to\mathbf{P}$ is the projection, then the fiber
$\pi_{s}^{-1}(Q)$ is the closed set of points $\mathbf{a}\in\check{\mathbb{P}_{k}^{n}}$
such that either $\dim\mathrm{Sing}(V\cap H_{\mathbf{a}})\geq\dim\mathrm{Sing}V$
or $\dim(V\cap H_{\mathbf{a}})=\dim V$. 

Now let $F_{1},\ldots,F_{r}$ be forms as in the hypothesis, and let
$q$ be a prime such that (\ref{eq:dimSing}) is satisfied for $Q\in\mathbf{P}$
representing the tuple of (mod $q$)-reductions $((F_{1})_{q},\ldots,(F_{r})_{q})$.
Then $\pi_{s}^{-1}(Q)$ is defined in $\mathbb{P}_{k}^{n}$, where
$k=\kappa(Q)=\mathbb{F}_{q}$, by the specializations $\left.H_{i}^{s}\right|_{Q}$
of the multihomogeneous forms $H_{i}^{s}$. Applying Lemma \ref{lem:Bertini.}
we get that $\pi_{s}^{-1}(Q)\times\Spec\bar{k}$ is a proper closed
subset of $\mathbb{P}_{\bar{k}}^{n}$ (where $\bar{k}$ is an algebraic
closure of $k$). Therefore one of the forms $\left.H_{i}^{s}\right|_{Q}\in k[\xi_{0},\ldots,\xi_{n}]$
must be non-zero, so the form \[
\Phi_{q}(\xi_{0},\ldots,\xi_{n})=\left.H_{i}^{s}\right|_{Q}(\xi_{0},\ldots,\xi_{n})\]
has the desired properties. 

The last assertion of the lemma follows from the easy observation
that a polynomial of degree at most $q$ cannot vanish at every point
of $\mathbb{P}_{\mathbb{F}_{q}}^{n}$.
\end{proof}
The following lemma explores the new geometry arising from the Weyl
differencing in Section \ref{sec:Proof-of-the-Main-Result}. For a
polynomial $f(X_{1},\ldots,X_{n})$ we denote by $\nabla f$ the gradient
$\left(\frac{\partial f}{\partial X_{1}},\ldots,\frac{\partial f}{\partial X_{n}}\right)^{t}$
and by $\nabla^{2}f$ the Hessian matrix $\left(\frac{\partial f}{\partial X_{i}\partial X_{j}}\right)_{1\leq i,j\leq n}$.

\begin{lem}
\label{lem:HBlemma2}Let $G_{1},\ldots,G_{r}$ be homogeneous polynomials
in $\mathbb{Z}[X_{1},\ldots,X_{n}]$ of degrees $d_{1},\ldots,d_{r}$,
and let \[
V=\Proj\mathbb{Z}[X_{1},\ldots,X_{n}]/(G_{1},\ldots,G_{r}).\]
Let $q$ be a prime such that $q\nmid d_{i}$ for all $i=1,\ldots,r$
and suppose that $V_{q}$ is a non-singular complete intersection
of codimension $r$ in $\mathbb{P}_{\mathbb{F}_{q}}^{n-1}$.
\end{lem}
\begin{itemize}
\item [(i)]\label{enu:dimS}Let\begin{multline*}
S=\left\{ (\mathbf{x},\mathbf{y})\in\mathbb{P}_{\mathbb{F}_{q}}^{n-1}\times\mathbb{P}_{\mathbb{F}_{q}}^{n-1};\ \mathbf{y}\cdot\nabla G_{i}(\mathbf{x})=0,\ i=1,\ldots,r,\right.\\
\left.\mathrm{rank}\left(\mathbf{y}\cdot\nabla^{2}G_{i}(\mathbf{x})\right)_{1\leq i\leq r}<r\right\} .\end{multline*}
Then $\dim S\leq n-2.$
\item [(ii)]\label{enu:dimT_s}For $\mathbf{y}\in\mathbb{P}_{\mathbb{F}_{q}}^{n-1},$
let\begin{multline*}
S_{\mathbf{y}}=\left\{ \mathbf{x}\in\mathbb{P}_{\mathbb{F}_{q}}^{n-1};\ \mathbf{y}\cdot\nabla G_{i}(\mathbf{x})=0,\ i=1,\ldots,r,\right.\\
\left.\mathrm{rank}\left(\mathbf{y}\cdot\nabla^{2}G_{i}(\mathbf{x})\right)_{1\leq i\leq r}<r,\right\} .\end{multline*}
For $s=-1,0,1,\ldots,n-1,$ let $T_{s}=\left\{ \mathbf{y}\in\mathbb{P}_{\mathbb{F}_{q}}^{n-1};\ \dim S_{\mathbf{y}}\geq s\right\} .$
Then $T_{s}$ is Zariski closed and $\dim T_{s}\leq n-s-2.$ 
\item [(iii)]For each $s$, let $T_{s}^{(1)},T_{s}^{(2)},\ldots$ be the
irreducible components of $T_{s}$. Then\[
\sum_{j}\deg(T_{s}^{(j)})=O_{n,r,d_{1},\ldots,d_{r}}(1).\]

\end{itemize}
To prove Lemma \ref{lem:HBlemma2} we shall need the following lemma.

\begin{lem}
\label{lem:diagonal}Let $k$ be a field, and let $V$ be a closed
subscheme of $\mathbb{P}_{k}^{n}\times\mathbb{P}_{k}^{n}.$ Let $\Delta\subseteq\mathbb{P}^{n}\times\mathbb{P}^{n}$
be the diagonal, $\Delta=\left\{ (\mathbf{x},\mathbf{x});\ \mathbf{x}\in\mathbb{P}_{k}^{n}\right\} $.
If $\dim V\geq n,$ then $V\cap\Delta\ne\emptyset.$
\begin{proof}
Consider the rational map \[
f:\mathbb{P}^{2n+1}\dashrightarrow\mathbb{P}^{n}\times\mathbb{P}^{n}\]
 given by \[
\left(X_{0}:\ldots:X_{2n+1}\right)\mapsto\left((X_{0}:\ldots:X_{n}),(X_{n+1}:\ldots:X_{2n+1})\right).\]
 Its domain of definition is the Zariski open set $U:=\mathbb{P}^{2n+1}\setminus(L\cup M),$
where $L=\{ X_{0}=\ldots=X_{n}=0\}$ and $M=\{ X_{n+1}=\ldots=X_{2n+1}=0\}$.
Moreover, let $\hat{\Delta}$ be the variety in $\mathbb{P}^{2n+1}$
defined by $X_{0}=X_{n+1},\ldots,X_{n}=X_{2n+1}.$ Then $f$ is an
isomorphism between $\hat{\Delta}$ and $\Delta.$ Let $\hat{V}$
be the Zariski closure in $\mathbb{P}^{2n+1}$ of $f^{-1}(V).$ Then
\[
\dim\hat{V}=\dim V+1\geq n+1,\]
so that\[
\mathrm{codim}\hat{\Delta}+\mathrm{codim}\hat{V}\leq2n+1.\]
Thus, by the Projective Dimension Theorem \cite[Ex. 3.3.4]{Liu},
$\hat{\Delta}\cap\hat{V}$ is nonempty. But a point $P$ in this intersection
automatically lies in $U$, since $\hat{\Delta}\cap(L\cup M)$ is
empty, and we get a point $f(P)$ in $\Delta\cap V.$ 
\end{proof}
\end{lem}
\begin{proof}[Proof of Lemma \ref{lem:HBlemma2}]

(i) Assume that $\dim S\geq n-1.$ According to Lemma \ref{lem:diagonal},
we then must have $S\cap\Delta\ne\emptyset.$ Thus, suppose $(\mathbf{x},\mathbf{x})\in S\cap\Delta.$
By the definition of $S,$ we then have \[
\begin{cases}
\mathbf{x}\cdot\nabla G_{i}(\mathbf{x})=0,\ i=1,\ldots,r\\
\mathrm{rank}\left(\mathbf{x}\cdot\nabla^{2}G_{i}(\mathbf{x})\right)_{1\leq i\leq r}<r.\end{cases}\]
But $\mathbf{x}\cdot\nabla^{2}G_{i}(\mathbf{x})=\nabla(\mathbf{x}\cdot\nabla G_{i}(\mathbf{x}))$,
so by Euler's identity we have (since $q$ does not divide any of
the degrees of the $G_{i}$)\[
\begin{cases}
G_{i}(\mathbf{x})=0,\ i=1,\ldots,r\\
\mathrm{rank}\left(\nabla G_{i}(\mathbf{x})\right)_{1\leq i\leq r}<r.\end{cases}\]
Therefore, by the Jacobian Criterion, $\mathbf{x}$ is a singular
point of $V,$ in contradiction with the hypothesis.

(ii) Let $\pi:S\rightarrow\mathbb{P}^{n-1}$ be the projection onto
the second coordinate, $(\mathbf{x},\mathbf{y})\mapsto\mathbf{y}.$
Then $S_{\mathbf{y}}=\pi^{-1}(\mathbf{y})\times\{\mathbf{y}\}.$ The
fact that $T_{s}$ is closed follows from Chevalley's semicontinuity
theorem \cite[Cor 13.1.5]{EGAIV(III)}. Now let $S_{s}=S\cap\left(\mathbb{P}^{n-1}\times T_{s}\right)$
for each $s=-1,\ldots,n-1.$ Since $S_{s}$ is the disjoint union
of fibres \[
S_{s}=\bigcup_{\mathbf{y}\in T_{s}}\pi^{-1}(\mathbf{y}),\]
we have, by (i)\[
\dim T_{s}+s\leq\dim S_{s}\leq\dim S\leq n-2,\]
whence $\dim T_{s}\leq n-s-2.$ 

(iii) As in Lemma \ref{lem:Phi}, we shall let $\mathbb{P}_{i}$ be
the projective spaces parametrizing hypersurfaces of degree $d_{i}$
in $\mathbb{P}_{\mathbb{Z}}^{n}$, and put $\mathbf{P}=\mathbb{P}_{1}\times\ldots\times\mathbb{P}_{r}$.
Now, let

\begin{multline*}
\mathcal{S}=\left\{ (G_{1},\ldots,G_{r},\mathbf{x},\mathbf{y})\in\mathbf{P}\times\mathbb{P}_{\mathbb{Z}}^{n-1}\times\mathbb{P}_{\mathbb{Z}}^{n-1};\ \mathbf{y}\cdot\nabla G_{i}(\mathbf{x})=0,\ i=1,\ldots,r,\right.\\
\left.\mathrm{rank}\left(\mathbf{y}\cdot\nabla^{2}G_{i}(\mathbf{x})\right)_{1\leq i\leq r}<r,\right\} .\end{multline*}
Let $\tilde{\pi}:\mathcal{S}\to\mathbf{P}\times\mathbb{P}_{\mathbb{Z}}^{n-1}$
be the projection $(G_{1},\ldots,G_{r},\mathbf{x},\mathbf{y})\mapsto(G_{1},\ldots,G_{r},\mathbf{y})$,
and define for each $s$\[
\mathcal{T}_{s}=\left\{ \mathcal{P}=(G_{1},\ldots,G_{r},\mathbf{y});\dim\tilde{\pi}^{-1}(\mathcal{P})\geq s\right\} .\]
Then $\mathcal{T}_{s}$ is closed by Chevalley's theorem, so it is
defined in $\mathbf{P}\times\mathbb{P}_{\mathbb{Z}}^{n-1}$ by multihomogeneous
polynomials $H_{1},\ldots,H_{t}$ where $t=O_{n,r,d_{1},\ldots,d_{r}}(1)$.
Now we fix polynomials $G_{1},\ldots,G_{r}$ and a prime $q$. The
set $T_{s}$ is then defined in $\mathbb{P}_{\mathbb{F}_{q}}^{n-1}$by
$\left.H_{1}\right|_{G_{1},\ldots,G_{r}},\ldots,\left.H_{t}\right|_{G_{1},\ldots,G_{r}}$.
Now by Bézout's Theorem \cite[Ex. 8.4.6]{Fulton} we have\[
\sum_{j}\deg(T_{s}^{(j)})\leq\prod_{i}\deg(H_{i})\ll_{n,r,d_{1},\ldots,d_{r}}1.\]

\end{proof}

\section{\label{sec:modq}Points on Complete Intersections over $\mathbb{F}_{q}$}

The following result is well-known and trivial, but we include a proof
for the sake of completeness.

\begin{lem}
\label{lem:trivial} Let $X=\Spec\mathbb{F}_{q}[X_{1},\ldots,X_{n}]/(f_{1},\ldots,f_{\rho})$
be a closed subscheme of $\mathbb{A}_{\mathbb{F}_{q}}^{n}$, and let
$d=\max_{i}(\deg f_{i})$. Let $B\geq1$. Then, for any box $\sB=\left[a_{1}-b_{1},a_{1}+b_{1}\right]\times\ldots\times\left[a_{n}-b_{n},a_{n}+b_{n}\right]$,
with $|b_{i}|\leq B$, containing at most one representative of each
congruence class modulo $q$, we have\[
N(X,\sB,q)\ll_{n,\rho,d}B^{\dim X}.\]

\end{lem}
\begin{proof}
We identify $\mathbb{A}_{\mathbb{F}_{q}}^{n}$ with the open subset
$\{ X_{0}\neq0\}$ of $\mathbb{P}_{\mathbb{F}_{q}}^{n}$ and consider
the scheme-theoretic closure $Y$ of $X$ in $\mathbb{P}_{\mathbb{F}_{q}}^{n}$
defined by the homogenizations $F_{1},\ldots,F_{\rho}$ of $f_{1},\ldots,f_{\rho}$.
Then the sum $D_{X}$ of the degrees of the irreducible components
of $Y$ is at most $d^{\rho}$ by Bézout's Theorem \cite[Ex. 8.4.6]{Fulton}.
Thus it suffices to show that $N(X,\sB,q)\ll_{n,D_{X}}B^{\dim X}$
for every closed subscheme $X$. We prove this by induction over $\nu=\dim X$.
If $\nu=0$, then $\# X(\mathbb{F}_{q})\leq D_{X}$, so we are done.
Thus, suppose that $\nu\geq1$. Since $X$ has at most $D_{X}$ irreducible
components, it is enough to prove that $N(X',\sB,q)\ll_{n,D_{X}}B^{\nu}$
for an arbitrary irreducible component $X'$ of $X$. For some $i\in\{1,\ldots,n\}$,
all the hyperplanes $H_{a}:$$x_{i}=a$, where $a$ ranges over $\mathbb{F}_{q}$,
intersect $X'$ properly. Since $D_{X\cap H_{a}}\leq D_{X}$, the
induction hypothesis yields that $N(X'\cap H_{a},\sB,q)\ll_{n,D_{X}}B^{\nu-1}$
for each $a\in\mathbb{F}_{q}$ . Since we only need to consider at
most $2B$ values of $a$, we get\[
N(X',\sB,q)=\sum_{a}N(X'\cap H_{a},\sB,q)\leq2B\cdot O_{n,D_{X}}(B^{\nu-1})\ll_{n,D_{X}}B^{\nu},\]
as desired.
\end{proof}
Delignes work on the Weil Concectures \cite{Deligne} yields a sharp
asymptotic formula for the number of $\mathbb{F}_{q}$-points on a
non-singular projective complete intersection. In the paper by Hooley
\cite{Hooley} (with an appendix by Katz) an extension to the singular
case is proven. The following lemma is an affine reformulation of
Hooley's result.

\begin{lem}
\label{lem:Hooley-Deligne}Let $Y$ be a closed subscheme of $\mathbb{P}_{\mathbb{F}_{q}}^{n}$
that is a complete intersection of codimension $r\leq n$ and multidegree
$(d_{1},\ldots,d_{r})$. Let $Z=Y\cap\{ x_{0}=0\}$ and suppose that
$\dim Z=\dim Y-1$. Put $X=Y\setminus Z$ and $s=\dim\mathrm{Sing}Z$.
Then we have\[
\# X(\mathbb{F}_{q})=q^{n-r}+O_{n,d_{1},\ldots,d_{r}}(q^{(n-r+2+s)/2}).\]

\end{lem}
\begin{proof}
In case $n=r$ the lemma is a trivial consequence of Bézout's Theorem.
We may thus assume that $n>r$. By \cite[Appendix, Thm. 1]{Hooley}
we have\[
\# Z(\mathbb{F}_{q})=1+q+\ldots+q^{n-r-1}+O(q^{(n-r+s)/2}).\]
However, $s\geq\dim\Sing Y-1$ by Remark \ref{rem:katzlemma3}, so
by the same theorem we get\[
\# Y(\mathbb{F}_{q})=1+q+\ldots+q^{n-r}+O(q^{(n-r+2+s)/2}).\]
Subtracting these two equations, we get\[
\# X(\mathbb{F}_{q})=q^{n-r}+O(q^{(n-r+2+s)/2}),\]
as stated.
\end{proof}
The following result is a generalization of Theorem 3 in \cite{Heath-Brown}.
However, even in the case of a hypersurface we get a slightly sharper
estimate. The reason for this is the use of estimates by Katz \cite{Katz}
for {}``singular'' exponential sums. A similar application of those
results are found in a paper by Luo \cite{Luo}.

\subsubsection*{Notation}

For an element $\mathbf{x}=(x_{1},\ldots,x_{n})$ in $\mathbb{Z}^{n}$
we let $\mathbf{x}_{q}=(x_{1}+q\mathbb{Z},\ldots,x_{n}+q\mathbb{Z})\in\mathbb{F}_{q}^{n}$. 

\begin{thm}
\label{thm:H-BThm3}Let $W:\mathbb{R}^{n}\to\mathbb{R}$ be an infinitely
differentiable function, supported in a cube of side $2L$. Let $q$
be a prime and $B$ a real number with $1\leq B\ll_{L}q$. Let\[
X=\Spec\mathbb{Z}[X_{1},\ldots,X_{n}]/(f_{1},\ldots,f_{r}),\]
where the leading forms $F_{1},\ldots,F_{r}$ of $f_{1},\ldots,f_{r}$
are of degree at least $2$, and let \[
Z_{q}=\Proj\mathbb{Z}[X_{1},\ldots,X_{n}]/(q,F_{1},\ldots,F_{r}).\]
Assume that $\dim Z_{q}=n-1-r.$ Let $s=\dim\Sing Z_{q}$ and $d=\max_{i}(\deg F_{i})$.
Define a weighted counting function\[
N_{W}(X,B,q)=\sum_{\substack{\mathbf{x}\in\mathbb{Z}^{n}\\
\mathbf{x}_{q}\in X_{q}}
}W\left(\frac{1}{B}\mathbf{x}\right).\]
Then we have\begin{equation}
\begin{aligned}N_{W}(X,B,q) & =q^{-r}N_{W}(\mathbb{A}^{n},B,q)\\
 & +O_{n,d,L}\left(D_{2n}B^{s+1}q^{(n-r-s-2)/2}(B+q^{1/2})\right),\end{aligned}
\label{eq:H-BThm3}\end{equation}
where, for each natural number $k$, $D_{k}$ is the maximum over
$\mathbb{R}^{n}$ of all partial derivatives of $W$ of order $k$.
\end{thm}
\begin{proof}
We begin with some preparatory considerations, to justify the use
of Lemma \ref{lem:Hooley-Deligne} later in the proof. Let \[
Y_{q}=\Proj\mathbb{Z}[X_{0},\ldots,X_{n}]/(q,G_{1},\ldots,G_{r}),\]
where $G_{i}(X_{0},\ldots,X_{n})=X_{0}^{d_{i}}f_{i}(X_{1}/X_{0},\ldots,X_{n}/X_{0})$
for $i=1,\ldots,n$. Then $Z_{q}=Y_{q}\cap\left\{ X_{0}=0\right\} $
and $X_{q}=Y_{q}\setminus Z_{q}$. Moreover, since $\dim Z_{q}=n-1-r$
we must have $\dim Y_{q}=n-r$.

We shall follow the approach of Heath-Brown \cite{Heath-Brown} and
use induction with respect to $s$, starting with the case when $Z_{q}$
is non-singular, that is, when $s=-1$. In case $n-r\geq2$ we shall
use Katz' results. We begin, however, with two trivial cases. Suppose
firstly that $n-r=1$. Then \[
N_{W}(X,B,q)\ll_{n,L}D_{0}N(X,B,q)\ll_{n,d}D_{0}B\]
 by Lemma \ref{lem:trivial}, and \[
q^{-r}N_{W}(\mathbb{A}^{n},B,q)\ll_{n,L}D_{0}q^{-n+1}B^{n}\ll_{n,L}D_{0}B,\]
so \begin{align*}
N_{W}(X,B,q)-q^{-r}N_{W}(\mathbb{A}^{n},B,q) & \ll_{n,d,L}D_{2n}(B+q^{1/2})\end{align*}
as required for (\ref{eq:H-BThm3}). Next, suppose that $n-r=0$.
Also in this case the formula (\ref{eq:H-BThm3}) holds, since $N_{W}(X,B,q)\ll_{n,d,L}D_{0}$
and $q^{-r}N_{W}(\mathbb{A}^{n},B,q)\ll_{n,L}D_{0}q^{-n}B^{n}\ll_{n,L}D_{0}$,
whereas the error term required for (\ref{eq:H-BThm3}) is $D_{2n}(Bq^{-1/2}+1)$.

From now on, we assume that $n-r\geq2$. By the Poisson Summation
Formula we have\begin{align*}
N_{W}(X,B,q)= & \sum_{\mathbf{z}\in X_{q}}\sum_{\mathbf{u}\in\mathbb{Z}^{n}}W\left(\frac{1}{B}\mathbf{(z}+q\mathbf{u})\right)\\
= & \sum_{\mathbf{z}\in X_{q}}\left(\frac{B}{q}\right)^{n}\sum_{\mathbf{a}\in\mathbb{Z}^{n}}e_{q}(\mathbf{a}\cdot\mathbf{z})\hat{W}\left(\frac{B}{q}\mathbf{a}\right)\\
= & \left(\frac{B}{q}\right)^{n}\sum_{\mathbf{a}\in\mathbb{Z}^{n}}\hat{W}\left(\frac{B}{q}\mathbf{a}\right)\Sigma_{q}(\mathbf{a}),\end{align*}
where \[
\Sigma_{q}(\mathbf{a})=\sum_{\mathbf{z}\in X_{q}}e_{q}(\mathbf{a\cdot z}),\]
a sum which we shall now investigate. In case $\mathbf{a}\equiv\mathbf{0}\pmod q$,
we can use Lemma \ref{lem:Hooley-Deligne} to conclude that we have
\[
\Sigma_{q}(\mathbf{a})=\# X_{q}(\mathbb{F}_{q})=q^{n-r}+O_{n,d}(q^{(n-r+1)/2}).\]
Next we consider $\Sigma_{q}(\mathbf{a})$ for $\mathbf{a}\not\equiv\mathbf{0}\pmod q$.
Since $Z_{q}$ is a projective complete intersection of dimension
at least $1$, it is geometrically connected. Being non-singular,
it is thus geometrically integral. The hypothesis that $\deg F_{i}\geq2$
for all $i$ now implies that for each $\mathbf{a}\in\mathbb{F}_{q}^{n}\setminus\left\{ \mathbf{0}\right\} $
we have $\dim(Z_{q}\cap H_{\mathbf{a}})=n-r-2$, where $H_{\mathbf{a}}$
is the hyperplane defined by $\mathbf{a}\cdot\mathbf{x}=0$. Then,
by Theorems 23 and 24 in \cite{Katz}, we have\[
\Sigma_{q}(\mathbf{a})\ll q^{(n-r+1+\delta(\mathbf{a}))/2},\]
where $\delta(\mathbf{a})=\dim\Sing(Z_{q}\cap H_{\mathbf{a}}).$ Thus
we get\begin{equation}
\begin{aligned}N_{W}(X,B,q) & =\left(\frac{B}{q}\right)^{n}\left(\sum_{q\mid\mathbf{a}}\hat{W}\left(\frac{B}{q}\mathbf{a}\right)\left(q^{n-r}+O_{n,d}\left(q^{(n-r+1)/2}\right)\right)\right)\\
 & +O\left(\left(\frac{B}{q}\right)^{n}\sum_{\mathbf{a}\in\mathbb{Z}^{n}}\left|\hat{W}\left(\frac{B}{q}\mathbf{a}\right)\right|q^{(n-r+1+\delta(\mathbf{a}))/2}\right).\end{aligned}
\label{eq:W1}\end{equation}
The first term here equals\begin{equation}
\begin{gathered}\left(\frac{B}{q}\right)^{n}q^{n-r}\sum_{\mathbf{v}\in\mathbb{Z}^{n}}\hat{W}\left(B\mathbf{v}\right)+O_{n,d}\left(\left(\frac{B}{q}\right)^{n}q^{(n-r+1)/2}\sum_{\mathbf{v}\in\mathbb{Z}^{n}}\hat{W}\left(B\mathbf{v}\right)\right)\\
=q^{-r}N_{W}(\mathbb{A}^{n},B,q)+O_{n,d,L}\left(B^{n}q^{-(n+r-1)/2}\right),\end{gathered}
\label{eq:term1}\end{equation}
by the Poisson formula in the reverse direction and since $N_{W}(\mathbb{A}^{n},B,q)=O_{n,d,L}(B^{n}).$
In order to estimate the second term in (\ref{eq:W1}) we write\[
\sum_{\mathbf{a}\in\mathbb{Z}^{n}}\left|\hat{W}\left(\frac{B}{q}\mathbf{a}\right)\right|q^{(n-r+1+\delta(\mathbf{a}))/2}=\Sigma_{1}+\Sigma_{2},\]
where \begin{gather*}
\Sigma_{1}=\sum_{|\mathbf{a}|\leq q/2}\left|\hat{W}\left(\frac{B}{q}\mathbf{a}\right)\right|q^{(n-r+1+\delta(\mathbf{a}))/2}\text{ and}\\
\Sigma_{2}=\sum_{|\mathbf{a}|>q/2}\left|\hat{W}\left(\frac{B}{q}\mathbf{a}\right)\right|q^{(n-r+1+\delta(\mathbf{a}))/2}.\end{gather*}
It follows from a result of Zak (see \cite[Appendix, Thm. 2]{Hooley})
that $\delta(\mathbf{a})=-1$ or $0$ for all $\mathbf{a}$. By Lemma
\ref{lem:Phi}, all $\mathbf{a}$ for which $\delta(\mathbf{a})=0$
satisfy $\Phi(\mathbf{a})\equiv0\pmod q$ for a non-zero polynomial
$\Phi(\xi_{1},\ldots,\xi_{n})$ with integer coefficients, whose degree
is $O_{n,d}(1)$. Thus, let us split $\Sigma_{1}$ into two sums\[
\Sigma_{1}=\sum_{\substack{|\mathbf{a}|\leq q/2\\
\Phi(\mathbf{a})\equiv0(q)}
}\left|\hat{W}\left(\frac{B}{q}\mathbf{a}\right)\right|q^{(n-r+1)/2}+\sum_{\substack{|\mathbf{a}|\leq q/2\\
\Phi(\mathbf{a})\not\equiv0(q)}
}\left|\hat{W}\left(\frac{B}{q}\mathbf{a}\right)\right|q^{(n-r)/2}\]
and denote the first by $\Sigma_{11}$ and the second by $\Sigma_{12}$.
We observe that, since the infinitely differentiable function $W$
has compact support, we have an estimate $\left|\hat{W}(\mathbf{t})\right|\ll_{n,L}D_{k}\left|\mathbf{t}\right|^{-k}$
for $\left|\mathbf{t}\right|\geq1$ and any $k\geq0,$ and moreover
$D_{k}\ll_{n,L}D_{k+1}$ for every $k$. In particular, for any $t\in\mathbb{R}^{n}$
we have the estimate \begin{equation}
\left|\hat{W}(\mathbf{t})\right|\ll_{n,L}D_{k}\min(1,\left|\mathbf{t}\right|^{-k}),\ k\geq0\label{eq:Fourier}\end{equation}
Thus we get\begin{align*}
\sum_{\substack{\mathbf{\left|a\right|\leq}q/2\\
\Phi(\mathbf{a})\equiv0(q)}
}\left|\hat{W}\left(\frac{B}{q}\mathbf{a}\right)\right| & \ll_{n,L}D_{2n}\sum_{\substack{\mathbf{\left|a\right|\leq}q/2\\
\Phi(\mathbf{a})\equiv0(q)}
}\min\left(1,\left|\frac{B}{q}\mathbf{a}\right|^{-2n}\right).\end{align*}
Without loss of generality we can assume that $\xi_{n}$ occurs in
the polynomial $\Phi(\xi_{1},\ldots,\xi_{n})$. Then, for each fixed
determination of $a_{1},\ldots,a_{n}$, there are $O_{n,d}(1)$ values
for which $\Phi(a_{1},\ldots,a_{n})\equiv0\pmod q$, and we get \begin{align*}
\sum_{\substack{\mathbf{\left|a\right|\leq}q/2\\
\Phi(\mathbf{a})\equiv0(q)}
}\min\left(1,\left|\frac{B}{q}\mathbf{a}\right|^{-2n}\right) & =\sum_{\left|a_{1}\right|\leq q/2}\cdots\sum_{\left|a_{n-1}\right|\leq q/2}\sum_{\substack{\left|a_{n}\right|\leq q/2\\
\Phi(\mathbf{a})\equiv0(q)}
}\min\left(1,\left|\frac{B}{q}\mathbf{a}\right|^{-2n}\right)\\
 & \ll_{n,d}\prod_{i=1}^{n-1}\sum_{\left|a_{i}\right|\leq q/2}\min\left(1,\left|\frac{B}{q}a_{i}\right|^{-2}\right).\end{align*}
Now, for each $i=1,\ldots,n-1$ we have\[
\sum_{\left|a_{i}\right|\leq q/2}\min\left(1,\left|\frac{B}{q}a_{i}\right|^{-2}\right)=\sum_{\left|a_{i}\right|\leq q/B}1+\sum_{q/B<\left|a_{i}\right|\leq q/2}\left|\frac{B}{q}a_{i}\right|^{-2}\ll\frac{q}{B},\]
and we conclude that\[
\Sigma_{11}\ll_{n,d,L}D_{2n}\left(\frac{q}{B}\right)^{n-1}q^{(n-r+1)/2}.\]
Moreover, using (\ref{eq:Fourier}) and the fact that\begin{equation}
\sum_{\substack{\mathbf{u}\in\mathbb{Z}^{n}\\
|\mathbf{u}|>U}
}|\mathbf{u}|^{-(n+1)}\ll_{n}U^{-1}\label{eq:sumU}\end{equation}
 we have\begin{align*}
\mathrm{\Sigma_{12}} & \leq\sum_{\substack{\mathbf{\left|a\right|\leq}q/2}
}\left|\hat{W}\left(\frac{B}{q}\mathbf{a}\right)\right|q^{(n-r)/2}\\
 & \leq q^{(n-r)/2}\left(\sum_{\substack{\mathbf{\left|a\right|\leq}q/B}
}\left|\hat{W}\left(\frac{B}{q}\mathbf{a}\right)\right|+\sum_{\substack{q/B<\mathbf{\left|a\right|\leq}q/2}
}\left|\hat{W}\left(\frac{B}{q}\mathbf{a}\right)\right|\right)\\
 & \ll_{n,L}D_{n+1}\left(\frac{q}{B}\right)^{n}q^{(n-r)/2}.\end{align*}
We arrive at the estimate\begin{equation}
\Sigma_{1}\ll_{n,d,L}D_{2n}\left(\frac{q}{B}\right)^{n}q^{(n-r-1)/2}(B+q^{1/2}).\label{eq:Sigma1}\end{equation}
It turns out that $\Sigma_{2}$ does not contribute to the error term.
Indeed, using (\ref{eq:Fourier}) and (\ref{eq:sumU}) again we have\begin{align*}
\Sigma_{2} & \leq\sum_{|\mathbf{a}|>q/2}\left|\hat{W}\left(\frac{B}{q}\mathbf{a}\right)\right|q^{(n-r+1)/2}\ll_{n,L}D_{n+1}\left(\frac{q}{B}\right)^{n}q^{(n-r-1)/2},\end{align*}
which is dominated by the bound (\ref{eq:Sigma1}) for $\Sigma_{1}$.
Thus, inserting (\ref{eq:term1}) and (\ref{eq:Sigma1}) into the
formula (\ref{eq:W1}) yields\begin{align*}
N_{W}(X,B,q) & =q^{-r}N_{W}(\mathbb{A}^{n},B,q)+O_{n,d,L}\left(D_{2n}q^{(n-r-1)/2}(B+q^{1/2})\right),\end{align*}
as required for the case $s=-1$.

Suppose now that $Z_{q}$ is singular, so that $s\geq0$. Following
Heath-Brown \cite{Heath-Brown} we will count points on hyperplane
sections. We begin with remarking that it is enough to prove the theorem
for $q$ greater than some constant $q_{0}=q_{0}(n,d)$. Indeed, if
$q\ll_{n,d}1$, then $B\ll_{n,d,L}1$, so that trivially we have $N_{W}(X,B,q)-q^{-r}N_{W}(\mathbb{A}^{n},B,q)\ll_{n,d,L}1.$
Thus, using Lemma \ref{lem:Phi}, we can assume that it is possible
to find a primitive integer vector $\mathbf{b}$, with $\mathbf{b}\ll_{n,d}1$,
such that $\dim(Z_{q}\cap H_{\mathbf{b}})=n-r-2$ and $\dim\Sing((Z_{q}\cap H_{\mathbf{b}})_{q})=s-1$,
where $H_{\mathbf{b}}$ is the hyperplane in $\mathbb{P}^{n-1}$ defined
by $\mathbf{b}\cdot\mathbf{x}=0$. We can find a unimodular integer
matrix $M$, all of whose entries are $O_{n,d}(1)$ such that the
automorphism of $\mathbb{P}_{\mathbb{Z}}^{n-1}$ induced by $M$ maps
$H_{\mathbf{b}}$ onto the hyperplane $X_{n}=0$, which we identify
with $\mathbb{P}^{n-2}=\Proj\mathbb{Z}[X_{1},\ldots,X_{n-1}]$. Let
$\tilde{Z}_{q}$ be the image of $Z_{q}\cap H_{\mathbf{b}}$. Then\[
\tilde{Z}_{q}=\Proj\mathbb{Z}[X_{1},\ldots,X_{n-1}]/(q,G_{1},\ldots,G_{r})\]
where $G_{i}(X_{1},\ldots,X_{n-1})=F_{i}(M^{-1}(X_{1},\ldots,X_{n-1},0))$
for $i=1,\ldots,r$, and each $G_{i}$ is of the same degree as $F_{i}$.
Obviously we have $\dim\Sing\tilde{Z}_{q}=s-1$. Moreover,\[
N_{W}(X,B,q)=\sum_{\mathbf{x}_{q}\in\mathbb{X}_{q}}W\left(\frac{1}{B}\mathbf{x}\right)=\sum_{\mathbf{x}_{q}\in\tilde{X}_{q}}\tilde{W}\left(\frac{1}{B}\mathbf{x}\right),\]
where $\tilde{X}$ is the image of $X$ under the automorphism of
$\mathbb{A}^{n}$ induced by $M$ and where $\tilde{W}(\mathbf{t})=W(M^{-1}\mathbf{t})$.
Then $\tilde{W}$ is supported in a cube of side $L'\ll_{n,d}L$,
so we can write\begin{equation}
N_{W}(X,B,q)=\sum_{-BL'\leq c\leq BL'}\sum_{\substack{\mathbf{x}_{q}\in\tilde{X}_{q}\\
x_{n}=c}
}\tilde{W}\left(\frac{1}{B}\mathbf{x}\right).\label{eq:slice-sum}\end{equation}
For each $c\in\mathbb{Z}$, the intersection of $\tilde{X}$ with
the hyperplane $x_{n}=c$ is isomorphic to \[
\tilde{X_{c}}=\Spec\mathbb{Z}[X_{1},\ldots,X_{n-1}]/(g_{1}^{c},\ldots,g_{r}^{c})\]
where $g_{i}^{c}(X_{1},\ldots,X_{n-1})=f_{i}(X_{1},\ldots,X_{n-1},c)$
for $i=1,\ldots,r$. For each $c$ and $i$, the leading form of $g_{i}^{c}$
is $G_{i}$, so our induction assumption applies to $\tilde{X}_{c},\tilde{Z}_{q}$
and the new weight function $\tilde{W}_{c}$ on $\mathbb{R}^{n-1}$
defined by $\tilde{W}_{c}(\mathbf{t})=\tilde{W}(\mathbf{t},c)$. We
get\begin{gather*}
\sum_{\substack{\mathbf{x}_{q}\in\tilde{X}_{q}\\
x_{n}=c}
}\tilde{W}\left(\frac{1}{B}\mathbf{x}\right)=N_{\tilde{W}_{c}}(\tilde{X}_{c},B,q)\\
=q^{-r}N_{\tilde{W}_{c}}(\mathbb{A}^{n-1},B,q)+O_{n,d,L}\left(D_{2n}B^{s}q^{(n-r-s-2)/2}(B+q^{1/2})\right).\end{gather*}
We shall now add the contributions from all $c$ in the interval $[-BL',BL']$.
Observe that\begin{align*}
\sum_{-BL'\leq c\leq BL'}N_{\tilde{W}_{c}}(\mathbb{A}^{n-1},B,q) & =\sum_{-BL'\leq c\leq BL'}\sum_{\mathbf{y}\in\mathbb{Z}^{n-1}}\tilde{W}\left(\frac{1}{B}(\mathbf{y},c)\right)\\
 & =\sum_{\mathbf{x}\in\mathbb{Z}^{n}}W\left(\frac{1}{B}M^{-1}\mathbf{x}\right)=\sum_{\mathbf{x}'\in\mathbb{Z}^{n}}W\left(\frac{1}{B}\mathbf{x}'\right)\\
 & =N_{W}(\mathbb{A}^{n},B,q),\end{align*}
since $M$ is unimodular. Thus, summing according to (\ref{eq:slice-sum})
we deduce that\[
N_{W}(X,B,q)=q^{-r}N_{W}(\mathbb{A}^{n},B,q)+O_{n,d,L}\left(D_{2n}B^{s+1}q^{(n-r-s-2)/2}(B+q^{1/2})\right)\]
and the induction step is finished.
\end{proof}

\section{\label{sec:Proof-of-the-Main-Result}Proof of the Main Result}

The aim of this section is to prove Theorem \ref{thm:main}. Throughout
the proof, any implicit constant is allowed to depend only on $n$
and $d$, and we will omit the subscripts $n,d$ from the $O$- and
$\ll$-notation. 

\begin{note*}
It will suffice to prove the theorem under the somewhat weaker hypothesis
that $p<2B+1<q$, but with the additional assumption that $2B+1$
is a multiple of $p$. We will now prove that the general case follows
from this case. If $p$ and $q$ are given primes and $B$ is an arbitrary
real number such that $2p<2B+1<q-p$, then there are integers $B_{1}$
and $B_{2}$, with $B_{1}\leq B\leq B_{2}$, such that $2B_{1}+1$
and $2B_{2}+1$ are multiples of $p$ and $p<2B_{i}+1<q$ for $i=1,2.$
We have \begin{align*}
N(X,B,pq)-\frac{(2B+1)^{n}}{p^{r}q^{r}} & \leq N(X,B_{2},pq)-\frac{(2B+1)^{n}}{p^{r}q^{r}}\\
 & =N(X,B_{2},pq)-\frac{(2B_{2}+1)^{n}}{p^{r}q^{r}}+O(B^{n-1}p^{-r+1}q^{-r}),\end{align*}
and similarly\[
N(X,B,pq)-\frac{(2B+1)^{n}}{p^{r}q^{r}}\geq N(X,B_{1},pq)-\frac{(2B_{1}+1)^{n}}{p^{r}q^{r}}+O(B^{n-1}p^{-r+1}q^{-r}).\]
Thus, if we assume Theorem \ref{thm:main} to be true for $B_{1}$
and $B_{2}$, then we see that it must also hold for $B$, since $B_{1},B_{2}\asymp B$.
\end{note*}
From now on we assume that $2B+1$ is a multiple of $p$ between $p$
and $q$. To facilitate the notation we introduce the characteristic
function of the box $\mathsf{B}=[-B,B]^{n}\cap\mathbb{Z}^{n}$,\[
\chi_{\mathsf{B}}(\mathbf{x})=\begin{cases}
1 & \text{if }\max|x_{i}|\leq B,\\
0 & \text{otherwise.}\end{cases}\]
Then \[
N:=N(X,B,pq)=\sum_{\substack{\mathbf{x}\in\mathbb{Z}^{n}\\
pq\mid f_{i}(\mathbf{x})}
}\chi_{\mathsf{B}}(\mathbf{x})=\sum_{\substack{\mathbf{w}\in\mathbb{F}_{p}^{n}\\
p\mid f_{i}(\mathbf{w})}
}\sum_{\substack{\mathbf{x}\equiv\mathbf{w}(p)\\
q\mid f_{i}(\mathbf{x})}
}\chi_{\mathsf{B}}(\mathbf{x}).\]
The {}``expected value'' of the inner sum is\[
K:=p^{-n}q^{-r}(2B+1)^{n},\]
so let us write\[
N=\sum_{\substack{\mathbf{w}\in\mathbb{F}_{p}^{n}\\
p\mid f_{i}(\mathbf{w})}
}\left(\sum_{\substack{\mathbf{x}\equiv\mathbf{w}(p)\\
q\mid f_{i}(\mathbf{x})}
}\chi_{\mathsf{B}}(\mathbf{x})-K\right)+K\sum_{\substack{\mathbf{w}\in\mathbb{F}_{p}^{n}\\
p\mid f_{i}(\mathbf{w})}
}1.\]
If we denote the first of these two sums by $S$, then, using Lemma
\ref{lem:Hooley-Deligne}, we get\begin{equation}
\begin{aligned}N & =S+K\# X(\mathbb{F}_{p})=S+K\left(p^{n-r}+O(p^{(n-r+1)/2})\right)\\
 & =\frac{(2B+1)^{n}}{p^{r}q^{r}}+S+O(B^{n}p^{-(n+r-1)/2}q^{-r}).\end{aligned}
\label{eq:N}\end{equation}
Now we turn our attention to $S.$ By Cauchy's inequality\begin{eqnarray*}
S^{2} & \leq & \left(\sum_{\substack{\mathbf{w}\in\mathbb{F}_{p}^{n}\\
p\mid f_{i}(\mathbf{w})}
}1\right)\left(\sum_{\substack{\mathbf{w}\in\mathbb{F}_{p}^{n}\\
p\mid f_{i}(\mathbf{w})}
}\left(\sum_{\substack{\mathbf{x}\equiv\mathbf{w}(p)\\
q\mid f_{i}(\mathbf{x})}
}\chi_{\mathsf{B}}(\mathbf{x})-K\right)^{2}\right),\end{eqnarray*}
so that, if we denote the expression in the rightmost parentheses
by $\Sigma$, and apply Lemma \ref{lem:trivial}, we get\begin{equation}
S\ll p^{(n-r)/2}\Sigma^{1/2}.\label{eq:S}\end{equation}
We estimate $\Sigma$ by adding some extra (positive) terms:\begin{eqnarray*}
\Sigma & \leq & \sum_{\substack{\mathbf{w}\in\mathbb{F}_{p}^{n}}
}\sum_{\mathbf{a}\in\mathbb{F}_{q}^{r}}\left(\sum_{\substack{\mathbf{x}\equiv\mathbf{w}(p)\\
f_{i}(\mathbf{x})\equiv a_{i}(q)}
}\chi_{\mathsf{B}}(\mathbf{x})-K\right)^{2}\\
 & = & \sum_{\substack{\mathbf{w}\in\mathbb{F}_{p}^{n}}
}\sum_{\mathbf{a}\in\mathbb{F}_{q}^{r}}\left(\sum_{\substack{\mathbf{x}\equiv\mathbf{w}(p)\\
f_{i}(\mathbf{x})\equiv a_{i}(q)}
}\chi_{\mathsf{B}}(\mathbf{x})\right)^{2}-2K\sum_{\mathbf{x}\in\mathbb{Z}^{n}}\chi_{\mathsf{B}}(\mathbf{x})+p^{n}q^{r}K^{2}.\end{eqnarray*}
The middle term here is just $-2p^{n}q^{r}K^{2}$, so, denoting the
first sum by $\mathcal{Z}$ we get\begin{equation}
\Sigma\leq\mathcal{Z}-p^{n}q^{r}K^{2}.\label{eq:Sigma}\end{equation}
To analyze $\mathcal{Z}$, we write\[
\mathcal{Z}=\sum_{\substack{\mathbf{x}\in\mathbb{Z}^{n}}
}\chi_{\mathsf{B}}(\mathbf{x})\sum_{\substack{\mathbf{x}'\in\mathbb{Z}^{n}\\
\mathbf{x}'\equiv\mathbf{x}(p)\\
f_{i}(\mathbf{x}')\equiv f_{i}(\mathbf{x})(q)}
}\chi_{\mathsf{B}}(\mathbf{x}').\]
 We make the variable change $\mathbf{x}'=\mathbf{x}+p\mathbf{y}$
in the second sum, introducing the {}``differentiated'' polynomials\[
f_{i}^{\mathbf{y}}(\mathbf{x})=f_{i}(\mathbf{x}+p\mathbf{y})-f_{i}(\mathbf{x}).\]
 If $\mathsf{B}_{\mathbf{y}}$ denotes the new box $\sB\cap(\mathsf{B}-p\mathbf{y})=\left\{ \mathbf{x}\in\mathbb{Z}^{n};\mathbf{x}\in\sB,\mathbf{x}+p\mathbf{y}\in\mathsf{B}\right\} $,
we get\begin{eqnarray*}
\mathcal{Z} & = & \sum_{\substack{\mathbf{x}\in\mathbb{Z}^{n}}
}\chi_{\mathsf{B}}(\mathbf{x})\sum_{\substack{\mathbf{y}\in\mathbb{Z}^{n}\\
f_{i}^{\mathbf{y}}(\mathbf{x})\equiv0(q)}
}\chi_{\mathsf{B}}(\mathbf{x}+p\mathbf{y})\\
 & = & \sum_{\substack{\mathbf{y}\in\mathbb{Z}^{n}}
}\sum_{\substack{\mathbf{x}\in\mathbb{Z}^{n}\\
f_{i}^{\mathbf{y}}(\mathbf{x})\equiv0(q)}
}\chi_{\mathsf{B}_{\mathbf{y}}}(\mathbf{x}).\end{eqnarray*}
Let us define\[
\Delta(\mathbf{y})=\sum_{\substack{\mathbf{x}\in\mathbb{Z}^{n}\\
f_{i}^{\mathbf{y}}(\mathbf{x})\equiv0(q)}
}\chi_{\mathsf{B}_{\mathbf{y}}}(\mathbf{x})-q^{-r}\sum_{\substack{\mathbf{x}\in\mathbb{Z}^{n}}
}\chi_{\mathsf{B}_{\mathbf{y}}}(\mathbf{x}),\]
and write\[
\mathcal{Z}=\sum_{\mathbf{y}\in\mathbb{Z}^{n}}\Delta(\mathbf{y})+q^{-r}\sum_{\mathbf{y}\in\mathbb{Z}^{n}}\sum_{\substack{\mathbf{x}\in\mathbb{Z}^{n}}
}\chi_{\mathsf{B}_{\mathbf{y}}}(\mathbf{x}).\]
Now one sees that, since we are assuming $p\mid(2B+1)$,\begin{align*}
\sum_{\mathbf{y}\in\mathbb{Z}^{n}}\sum_{\substack{\mathbf{x}\in\mathbb{Z}^{n}}
}\chi_{\mathsf{B}_{\mathbf{y}}}(\mathbf{x}) & =\prod_{i=1}^{n}\left(\sum_{y_{i}\in\mathbb{Z}}\sum_{\substack{x_{i}\in\mathbb{Z}}
}\chi_{[-B,B]}(x_{i})\chi_{[-B-py_{i},B-py_{i}]}(x_{i})\right)\\
 & =\left(\frac{(2B+1)^{2}}{p}\right)^{n}=p^{n}q^{2r}K^{2}.\end{align*}
In other words, $\mathcal{Z}=\sum\Delta(\mathbf{y})+p^{n}q^{r}K^{2}$,
so we get by (\ref{eq:Sigma}) \begin{equation}
\Sigma\leq\sum_{\mathbf{y}\in\mathbb{Z}^{n}}\Delta(\mathbf{y}).\label{eq:Sigma2}\end{equation}
 Our task is now to estimate $\sum\Delta(\mathbf{y})$. To this end,
denote the leading forms of $f_{1}^{\mathbf{y}},\ldots,f_{r}^{\mathbf{y}}$
by $F_{1}^{\mathbf{y}},\ldots F_{r}^{\mathbf{y}}$ and let \begin{gather*}
X_{\mathbf{y}}=\Spec\mathbb{F}_{q}[x_{1},\ldots,x_{n}]/(f_{1}^{\mathbf{y}},\ldots,f_{r}^{\mathbf{y}}),\\
Z_{\mathbf{y}}=\Proj\mathbb{F}_{q}[x_{1},\ldots,x_{n}]/(F_{1}^{\mathbf{y}},\ldots,F_{r}^{\mathbf{y}}).\end{gather*}
Observe that for each $i=1,\ldots,r$ we have\[
F_{i}^{\mathbf{y}}=p\mathbf{y}\cdot\nabla F_{i},\]
unless the right hand side vanishes identically (mod $q$) in $\mathbf{x}$.
Due to the non-singularity of $Z$, this happens only if $\mathbf{y}\equiv0\pmod q$.
Indeed, if $\mathbf{y}\cdot\nabla F_{i}$ is identically zero for
some $i$, then, in the notation of Lemma \ref{lem:HBlemma2}, $S_{\mathbf{y}}=\mathbb{P}_{\mathbb{F}_{q}}^{n-1}$.
Thus $\mathbf{y}$ is a point on the affine cone over $T_{n-1}=\emptyset$. 

\begin{lem}
\label{lem:Delta(y)}\begin{align*}
\sum_{\mathbf{y}\in\mathbb{Z}^{n}}\Delta(\mathbf{y}) & \ll B^{n+1}p^{-n}q^{(n-r-1)/2}(\log q)^{n}+B^{n+1}p^{-r}q^{-1/2}(\log q)^{n}\\
 & \phantom{\ll}+B^{n}p^{-n}q^{(n-r)/2}(\log q)^{n}+B^{n}(\log q)^{n}.\end{align*}

\end{lem}
\begin{proof}
First, we note that $\Delta(\mathbf{y})=0$ for all $\mathbf{y}$
with $\left|\mathbf{y}\right|\geq(2B+1)/p.$ Thus, we only need to
sum over the set\[
\mathcal{B}=\left\{ \mathbf{y}\in\mathbb{Z}^{n};\left|\mathbf{y}\right|<(2B+1)/p\right\} .\]
Let us decompose this set into subsets: $\mathcal{B}=\mathcal{B}_{0}\cup\mathcal{B}_{1}\cup\ldots\cup\mathcal{B}_{r}$,
where\[
\mathcal{B}_{\sigma}=\left\{ \mathbf{y}\in\mathcal{B};\mathrm{codim}Z_{\mathbf{y}}=\sigma\right\} ,\quad\sigma=0,\ldots,r.\]
For $\mathbf{y}\in\mathcal{B}_{r},$ we can use Theorem 1 of the Appendix
\cite{Salberger} to get\[
\Delta(\mathbf{y})\ll_{n,d}B^{s(\mathbf{y})+1}q^{(n-r-s(\mathbf{y})-2)/2}(B+q^{1/2})(\log q)^{n},\]
where $s(\mathbf{y})=\dim\mathrm{Sing}(Z_{\mathbf{y}}).$ Next we
need to find out how often each value of \textbf{$s(\mathbf{y})$}
arises. We consult Lemma \ref{lem:HBlemma2}. Since $Z_{\mathbf{y}}$
is a complete intersection of codimension $r$, the Jacobian Criterion
implies that $\mathrm{Sing}(Z_{\mathbf{y}})=S_{\mathbf{y}}.$ Thus,
the set of all \textbf{$\mathbf{y}$} such that $s(\mathbf{y})=s$
is contained in the affine cone over the set $T_{s}$. By part (ii)
of Lemma \ref{lem:HBlemma2}, $T_{s}$ has projective dimension $n-s-2$,
so by part (iii) and Lemma \ref{lem:trivial}, we get\[
\#\left\{ \mathbf{y}\in\mathcal{B}_{r};s(\mathbf{y})=s\right\} \ll_{n,d}\left(\frac{B}{p}\right)^{n-s-1}.\]
Summing, we get\begin{multline*}
\sum_{\mathbf{y}\in\mathcal{B}_{r}}\Delta(\mathbf{y})\ll\sum_{s=-1}^{n-r-1}\left(\frac{B}{p}\right)^{n-s-1}B^{s+1}q^{(n-r-s-2)/2}(B+q^{1/2})(\log q)^{n}\\
\ll B^{n}(\log q)^{n}\left(Bp^{-n}q^{(n-r-1)/2}+p^{-n}q^{(n-r)/2}+Bp^{-r}q^{-1/2}+p^{-r}\right).\end{multline*}
It remains to consider the contribution from $\mathbf{y}\in\mathcal{B}_{\sigma},\ \sigma<r.$
We make a simple observation about the varieties $Z_{\mathbf{y}}$
originating from these values of $\mathbf{y}$: now the set $S_{\mathbf{y}}$
is very large.
\begin{lem}
\label{lem:Jacobian}Let $G_{1},\ldots,G_{r}$ be forms in the variables
$X_{1},\ldots,X_{n}.$ Let\[
V=\{ G_{1}=\ldots=G_{r}=0\}\subseteq\mathbb{P}^{n-1}\]
and let\[
W=\left\{ G_{1}=\ldots=G_{r}=0,\ \mathrm{rank}\left(\frac{\partial G_{i}}{\partial X_{j}}\right)<r\right\} .\]
Suppose that $\mathrm{codim}(V)=\sigma<r$. Then $W$ contains all
irreducible components of $V$ of dimension $n-1-\sigma$. In particular,
$\dim W=n-1-\sigma$. 
\end{lem}
\begin{proof}
Let $V'$ be an irreducible component of $V$ with $\dim V'=n-1-\sigma$.
Assume that there were a point $P\in V'$ such that $\mathrm{rank}\left(\frac{\partial G_{i}}{\partial X_{j}}\right)(P)=r.$
Then we would have \[
\dim T_{P}V'=n-1-r<n-1-\sigma=\dim V',\]
 a contradiction. Thus $V'\subseteq W.$
\end{proof}
We see that if $\mathbf{y}\in\mathcal{B}_{\sigma}$, then, by Lemma
\ref{lem:Jacobian}, $\dim S_{\mathbf{y}}=n-1-\sigma$. Recalling
that, in the notation of Lemma \ref{lem:HBlemma2}, $T_{n-1-\sigma}$
has dimension less than or equal to $\sigma-1$, we must have\[
\left|\mathcal{B}_{\sigma}\right|\ll\left(\frac{B}{p}\right)^{\sigma}.\]
Using Lemma \ref{lem:trivial} to get the trivial estimate $\Delta(\mathbf{y})\ll B^{n-\sigma}$
for $\mathbf{y}\in\mathcal{B}_{\sigma},$ we compute the contribution
from the $\mathcal{B}_{\sigma}$, $\sigma<r$:\[
\sum_{\sigma=0}^{r-1}\sum_{\mathbf{y}\in\mathcal{B}_{\sigma}}\Delta(\mathbf{y)}=\sum_{\sigma=0}^{r-1}\left(\frac{B}{p}\right)^{\sigma}B^{n-\sigma}=B^{n}\sum_{\sigma=0}^{r-1}p^{-\sigma}\ll B^{n}.\]
In sum, then,\begin{align*}
\sum_{\mathbf{y}\in\mathcal{B}}\Delta(\mathbf{y)} & =\sum_{\sigma=0}^{r}\sum_{\mathbf{y}\in\mathcal{B}_{\sigma}}\Delta(\mathbf{y)}\\
 & \ll B^{n+1}p^{-n}q^{(n-r-1)/2}(\log q)^{n}+B^{n+1}p^{-r}q^{-1/2}(\log q)^{n}\\
 & \phantom{\ll}+B^{n}p^{-n}q^{(n-r)/2}(\log q)^{n}+B^{n}(\log q)^{n},\end{align*}
and Lemma \ref{lem:Delta(y)} follows.
\end{proof}
Working our way back through the estimates (\ref{eq:Sigma2}), (\ref{eq:S})
and (\ref{eq:N}), we now arrive at

\begin{multline}
\begin{gathered}N=\frac{(2B+1)^{n}}{p^{r}q^{r}}+O\left(B^{(n+1)/2}p^{-r/2}q^{(n-r-1)/4}(\log q)^{n/2}\right.\\
+B^{(n+1)/2}p^{(n-2r)/2}q^{-1/4}(\log q)^{n/2}+B^{n/2}p^{-r/2}q^{(n-r)/4}(\log q)^{n/2}\\
\left.+B^{n/2}p^{(n-r)/2}(\log q)^{n/2}+B^{n}p^{-(n+r-1)/2}q^{-r}\right).\end{gathered}
\label{eq:arrive}\end{multline}
 This completes the proof of Theorem \ref{thm:main}.

We shall now prove Corollary \ref{thm:uniform}, where the modest
dependence upon $\left\Vert F_{i}\right\Vert $ is due to the following
lemma.

\begin{lem}
\label{lem:goodp}Let $X$ and $Z$ be defined as in Theorem \ref{thm:main},
and assume that $Z_{\mathbb{Q}}$ is non-singular of dimension $n-1-r$.
If $P\geq\left(\sum_{i=1}^{r}\log\left\Vert F_{i}\right\Vert \right)^{1+\delta}$,
then there is a prime $p\asymp_{\delta}P$ such that $Z_{p}$ is non-singular
of dimension $n-1-r$.
\end{lem}
\begin{proof}
As in the proof of Lemma \ref{lem:Phi}, let $\mathbf{P}=\mathbb{P}_{1}\times\ldots\times\mathbb{P}_{r}$,
where $\mathbb{P}_{i}$ is the projective space parametrizing all
hypersurfaces of degree $d_{i}$ in $\mathbb{P}_{\mathbb{Z}}^{n-1}$.
By a semicontinuity argument analogous to that in the proof of Lemma
\ref{lem:Phi}, the subset $U\subseteq\mathbf{P}$ defined by \[
(G_{1},\ldots,G_{r})\in U\Leftrightarrow V(G_{1},\ldots,G_{r})\text{ is non-singular of codimension }r,\]
is Zariski open, its complement thus being defined by multihomogeneous
polynomials $H_{1},\ldots,H_{t}$ in the coefficients of $G_{1},\ldots,G_{r}$.
Now by the hypotheses, for some $j$ we must have $H_{j}(F_{1},\ldots,F_{r})\neq0$.
We observe firstly that\[
\log\left|H_{j}(F_{1},\ldots,F_{r})\right|\ll_{n,d}\sum_{i=1}^{r}\log\left\Vert F_{i}\right\Vert .\]
Secondly, for an arbitrary positive number $A$ we have\[
\#\left\{ p>AP;p\mid H_{j}(F_{1},\ldots,F_{r})\right\} \ll\frac{\log\left|H_{j}(F_{1},\ldots,F_{r})\right|}{\log AP}.\]
Thus, if we choose $A$ large enough, there are fewer than \[
a:=\left[\sum_{i=1}^{r}\log\left\Vert F_{i}\right\Vert \right]\]
 such primes. Hence among the $a$ first prime numbers greater than
$AP$, there must be one prime $p$ such that $p\nmid H_{j}(F_{1},\ldots,F_{r})$.
By Chebyshev's Theorem it is possible to find an interval $[AP,c_{\delta}AP]$
that contains more than $P^{1/(1+\delta)}$ primes. Since $P\geq a^{1+\delta}$,
this interval must contain $p$.
\end{proof}
Now we are ready to prove Theorem \ref{thm:uniform}.

\begin{proof}
[Proof of Theorem \ref{thm:uniform}]Theorem \ref{thm:main} yields
in particular that\begin{multline*}
N(X,B,pq)\ll_{n,d}\left[\frac{B^{n}}{p^{r}q^{r}}+B^{(n+1)/2}p^{-r/2}q^{(n-r-1)/4}\right.\\
+B^{(n+1)/2}p^{(n-2r)/2}q^{-1/4}+B^{n/2}p^{-r/2}q^{(n-r)/4}+B^{n/2}p^{(n-r)/2}\\
\left.+B^{n}p^{-(n+r-1)/2}q^{-r}+B^{n-1}p^{-r+1}q^{-r}\vphantom{\frac{B^{n}}{p^{r}q^{r}}}\right](\log q)^{n/2}\end{multline*}
Thus we want to optimize the expression\begin{multline*}
\frac{B^{n}}{p^{r}q^{r}}+B^{(n+1)/2}p^{-r/2}q^{(n-r-1)/4}+B^{(n+1)/2}p^{(n-2r)/2}q^{-1/4}\\
+B^{n/2}p^{-r/2}q^{(n-r)/4}+B^{n/2}p^{(n-r)/2}+B^{n}p^{-(n+r-1)/2}q^{-r}+B^{n-1}p^{-r+1}q^{-r}\end{multline*}
by choosing appropriate $p$ and $q$. It turns out that\begin{align}
p & \asymp B^{1-\frac{5nr-r^{2}-5r}{n^{2}+4nr-n-r^{2}-r}}, & q & \asymp B^{2-\frac{2(4nr-r^{2})}{n^{2}+4nr-n-r^{2}-r}}.\label{eq:pq}\end{align}
would be an optimal choice. (Note that the last two terms in the expression
are dominated by the first term, so the optimization consists of trying
to get the first five terms to be of approximately equal order of
magnitude.) The restriction $n\geq4r+2$ ensures that (\ref{eq:pq})
is compatible with the requirement that $2p<2B+1<q-p$. The trouble
is now to make sure that the intervals specified in (\ref{eq:pq})
contain {}``good'' primes, that is, primes such that both $Z_{p}$
and $Z_{q}$ are non-singular of dimension $n-1-r$. 

For $B$ large enough, (\ref{eq:pq}) is a valid choice. Indeed, if
\begin{gather*}
B\geq\left(\sum_{i=1}^{r}\log\left\Vert F_{i}\right\Vert \right)^{e_{1}}\text{, where }\\
e_{1}=\left(1-\frac{5nr-r^{2}-5r}{n^{2}+4nr-n-r^{2}-r}\right)^{-1}\left(1+\frac{1}{2r}\right),\end{gather*}
then by Lemma \ref{lem:goodp} (with $\delta=(2r)^{-1}$) we can choose
$p$ and $q$, satisfying (\ref{eq:pq}), such that Theorem \ref{thm:main}
holds. For these $B$, and with $p$ and $q$ subject to (\ref{eq:pq}),
Theorem \ref{thm:main} implies that\[
N(X,B)\ll_{n,d}N(X,B,pq)\ll_{n,d}B^{n-3r+r^{2}\frac{13n-3r-5}{n^{2}+4nr-n-r^{2}-r}}(\log B)^{n/2}.\]
 For $B<\left(\sum_{i=1}^{r}\log\left\Vert F_{i}\right\Vert \right)^{e_{1}}$,
we use the trivial estimate \[
N(X,B)\ll_{n,d}B^{n-r}\]
 obtained by Lemma \ref{lem:trivial} to get\begin{gather*}
N(X,B)\ll_{n,d}B^{n-3r+r^{2}\frac{13n-3r-5}{n^{2}+4nr-n-r^{2}-r}}\left(\sum_{i=1}^{r}\log\left\Vert F_{i}\right\Vert \right)^{e_{2}}\text{, where}\\
e_{2}=e_{1}\left(2r-r^{2}\frac{13n-3r-5}{n^{2}+4nr-n-r^{2}-r}\right)\leq2r+1.\end{gather*}
This proves the theorem.
\end{proof}
\begin{rem*}
If we are content with just an upper bound for $N(X,B,pq)$ in Theorem
\ref{thm:main}, we can get rid of the factor $(\log q)^{n/2}$ and
thus prove a slightly sharpened version of Theorem \ref{thm:uniform},
without the factor $(\log B)^{n/2}$. This can be achieved by introducing
an infinitely differentiable weight function into the proof of Theorem
\ref{thm:main}, as in \cite{Heath-Brown}, and using Theorem \ref{thm:H-BThm3}
in the place of \cite[Thm. 1]{Salberger} . More precisely, if instead
of $N(X,B,pq)$ we consider the weighted counting function\[
N_{W}(X,B,pq)=\sum_{\substack{\mathbf{x}\in\mathbb{Z}^{n}\\
\mathbf{x}_{p}\in X_{p}\\
\mathbf{x}_{q}\in X_{q}}
}W\left(\frac{1}{2B}\mathbf{x}\right),\]
where $W$ is a non-negative, infinitely differentiable weight function
on $\mathbb{R}^{n}$ supported in $[-1,1]^{n}$, we can prove an asymptotic
formula for $N_{W}(X,B,pq)$ where the main term is \[
p^{-r}q^{-r}\sum_{\mathbf{x}\in\mathbb{Z}^{n}}W\left(\frac{1}{2B}\mathbf{x}\right).\]
The error term would then consist of the first four error terms of
Theorem \ref{thm:main} with the factor $(\log q)^{n/2}$ removed,
the fifth error term unchanged, and an additional term which is $o\left(p^{-r}q^{-r}B^{n}\right)$
and thus negligible for the application of Theorem \ref{thm:uniform}.
To prove this asymptotic formula one imitates the proof of Theorem
\ref{thm:main}, with $\chi_{\sB}(\mathbf{x})$ replaced by $W\left(\frac{1}{2B}\mathbf{x}\right)$
and $K$ by \[
K_{W}=p^{-n}q^{-r}\sum_{\mathbf{x}\in\mathbb{Z}^{n}}W\left(\frac{1}{2B}\mathbf{x}\right).\]
One is then led to estimate expressions\[
\Delta_{W}(\mathbf{y})=\sum_{\substack{\mathbf{x}\in\mathbb{Z}^{n}\\
f_{i}^{\mathbf{y}}(\mathbf{x})\equiv0(q)}
}W_{\mathbf{y}}(\mathbf{x})-q^{-r}\sum_{\mathbf{x}\in\mathbb{Z}^{n}}W_{\mathbf{y}}(\mathbf{x}),\]
where $W_{\mathbf{y}}(\mathbf{x})=W\left(\frac{1}{2B}\mathbf{x}\right)W\left(\frac{1}{2B}(\mathbf{x}+p\mathbf{y})\right)$.
At this point we invoke Theorem \ref{thm:H-BThm3}. Here the error
term, in contrast to the unweighted formula of Theorem 1 in the Appendix,
contains no factor $(\log q)^{n}$, whence the promised improvement
of the upper bound. The only main divergence from the proof of Theorem
\ref{thm:main} lies in the calculation of the sum $\sum_{\mathbf{y}\in\mathbb{Z}^{n}}\sum_{\mathbf{x}\in\mathbb{Z}^{n}}W_{\mathbf{y}}(\mathbf{x})$.
This can be done by means of Poisson summation (see \cite[p. 20]{Heath-Brown})
and gives rise to the additional error term mentioned above. 
\end{rem*}

\bibliographystyle{plain}
\bibliography{ratpoints}

\newpage
\appendix
\section*{Appendix}
\begin{center}
\small\textsc{Per Salberger}
\end{center}

The aim of this note is to count $\mathbb{F}_{q}$-points in boxes
on affine varieties. If $\mathbf{x}=(x_{1},\ldots,x_{n})\in\mathbb{Z}^{n}$
and $q$ is a prime, then we set $\mathbf{x}_{q}=(x_{1}+q\mathbb{Z},\ldots,x_{n}+q\mathbb{Z})\in\mathbb{F}_{q}^{n}$.
If $\mathsf{B}$ is a box in $\mathbb{R}^{n}$ and $W$ a closed subscheme
of $\mathbb{A}_{\mathbb{Z}}^{n}$, then we let\[
N(W,\sB,q)=\#\left\{ \mathbf{x}=(x_{1},\ldots,x_{n})\in\sB\cap\mathbb{Z}^{n}:\mathbf{x}_{q}\in W(\mathbb{F}_{q})\right\} .\]

\begin{applem}
\label{lem:appendix}Let $q$ be a prime and $\sB$ be a box in $\mathbb{R}^{n}$
such that each side has length at most $2B<q$. Let $\fs,l_{1},\ldots,l_{s+1}$
be polynomials in $\mathbb{Z}[x_{1},\ldots,x_{n}],$ $r+s+1\leq n$
such that the leading forms $F_{1},\ldots,F_{r}$ of $\fs$ are of
degree $\geq2$ and the leading forms $L_{1},\ldots,L_{s+1}$ of $l_{1},\ldots,l_{s+1}$are
of degree $1$. Let \begin{gather*}
X=\Spec\mathbb{Z}[x_{1},\ldots,x_{n}]/(f_{1},\ldots,f_{r},l_{1}\ldots,l_{s+1}),\\
\Lambda=\Spec\mathbb{Z}[x_{1},\ldots,x_{n}]/(l_{1}\ldots,l_{s+1})\text{ and}\\
Z=\Proj\mathbb{Z}[x_{1},\ldots,x_{n}]/(F_{1},\ldots,F_{r},L_{1},\ldots,L_{s+1}).\end{gather*}
Suppose that $Z_{q}=Z_{\mathbb{F}_{q}}$ is non-singular of codimension
$r+s+1$ in $\mathbb{P}_{\mathbb{F}_{q}}^{n-1}$. Then\[
N(X,\sB,q)=q^{-r}N(\Lambda,\sB,q)+O_{n,d}(q^{(n-r-s-2)/2}(B+q^{1/2})(\log q)^{n}),\]
where $d=\max_{i}\deg F_{i}$. 
\end{applem}
\begin{proof}
If $r+s+1=n$, then $\# X(\mathbb{F}_{q})\leq d^{n}$ by the theorem
of Bezout and hence $N(X,\sB,q)-q^{-r}N(\Lambda,\sB,q)\ll_{n,d}1\leq q^{(n-r-s-2)/2}(B+q^{1/2}).$
If $r+s+1=n-1$, then $N(X,\sB,q)=O_{n,d}(B)$ by Lemma 5 in \cite{marmon}
so that $N(X,\sB,q)-q^{-r}N(\Lambda,\sB,q)\ll_{n,d}B\leq q^{(n-r-s-2)/2}(B+q^{1/2}).$
We may thus assume that $r+s+1\leq n-2$. Then, $Z_{q}$ is geometrically
connected since it is a complete intersection of dimension $\geq1$
(see \cite[Ex. II.8.4(c)]{hartshorne}). It is thus geometrically
integral since it is non-singular. Therefore, by the homogeneous Nullstellensatz
we obtain that a linear form $\mathbf{a}\cdot\mathbf{x}=a_{1}x_{1}+\ldots+a_{n}x_{n}$,
$(a_{1},\ldots,a_{n})\in\mathbb{F}_{q}^{n}$ vanishes on $Z_{q}$
if and only if $\mathbf{a}\cdot\mathbf{x}$ belongs to the linear
$\mathbb{F}_{q}$-space V of linear forms in $(x_{1},\ldots,x_{n})$
generated by the reductions of $L_{1},\ldots,L_{s+1}\pmod q$. We
now follow the approach of \cite{luo}. Let $S_{1}(\mathbf{a})=\sum_{\mathbf{b}\in\sB\cap\mathbb{Z}^{n}}e_{q}(-\mathbf{a}\cdot\mathbf{b})$
and $S_{2}(\mathbf{a})=\sum_{\mathbf{x}\in X(\mathbb{F}_{q})}e_{q}(\mathbf{a}\cdot\mathbf{x})$
for $\mathbf{a}\in\mathbb{F}_{q}$. Then,\[
N(X,\sB,q)=q^{-n}\sum_{\mathbf{a}\in\mathbb{F}_{q}^n}S_{1}(\mathbf{a})S_{2}(\mathbf{a}).\]
Let $\Pi_{\mathbf{a}}=\Proj\mathbb{F}_{q}[x_{1},\ldots,x_{n}]/(a_{1}x_{1}+\ldots+a_{n}x_{n})$
for $\mathbf{a}=(a_{1},\ldots,a_{n})\in\mathbb{F}_{q}^{n}$. Then,\begin{gather*}
q^{-(s+1)}\sum_{\mathbf{a}\in V}S_{1}(\mathbf{a})S_{2}(\mathbf{a})=q^{-(s+1)}\sum_{\mathbf{a}\in V}\sum_{\mathbf{x}\in X(\mathbb{F}_{q})}\sum_{\mathbf{b}\in\sB\cap\mathbb{Z}^{n}}e_{q}(\mathbf{a}\cdot(\mathbf{x}-\mathbf{b}))\\
=\sum_{\mathbf{x}\in X(\mathbb{F}_{q})}\sum_{\mathbf{b}\in\sB\cap\mathbb{Z}^{n}}\prod_{i=1}^{s+1}\left(\frac{1}{q}\sum_{a\in\mathbb{F}_{q}}e_{q}(aL_{i}(\mathbf{x}-\mathbf{b}))\right)\\
=\#\left\{ (\mathbf{x},\mathbf{b})\in X(\mathbb{F}_{q})\times(\sB\cap\mathbb{Z}^{n}):L_{1}(\mathbf{x}-\mathbf{b})\equiv\ldots\equiv L_{s+1}(\mathbf{x}-\mathbf{b})\equiv0\pmod q\right\} \\
=\#\left\{ (\mathbf{x},\mathbf{b})\in X(\mathbb{F}_{q})\times(\sB\cap\mathbb{Z}^{n}):l_{1}(\mathbf{b})\equiv\ldots\equiv l_{s+1}(\mathbf{b})\equiv0\pmod q\right\} \\
=\# X(\mathbb{F}_{q})N(\Lambda,\sB,q).\end{gather*}
Here $\# X(\mathbb{F}_{q})=q^{n-r-s-1}+O_{n,d}(q^{(n-r-s)/2})$ by
Lemma 6 in \cite{marmon}. There is also a set of $n-s-1$ indices
$i(1),\ldots,i(n-s-1)\in\{1,\ldots,n\}$ such that any $\mathbf{b}=(b_{1},\ldots,b_{n})\in\sB\cap\mathbb{Z}^{n}$
with $\mathbf{b}_{q}\in\Lambda(\mathbb{F}_{q})$ is uniquely determined
by $(b_{i(1)},\ldots,b_{i(n-s-1)})$. Hence, $\# N(\Lambda,\sB,q)\ll_{n}B^{n-s-1}$.
We have thus shown that \begin{align*}
q^{-n}\sum_{\mathbf{a}\in V}S_{1}(\mathbf{a})S_{2}(\mathbf{a}) & =q^{-(n-s-1)}\# X(\mathbb{F}_{q})N(\Lambda,\sB,q)\\
 & =q^{-r}N(\Lambda,\sB,q)+O_{n,d}(q^{-(n-s-1)+(n-r-s)/2}B^{n-s-1}).\end{align*}
As $q^{-(n-s-1)+(n-r-s)/2}B^{n-s-1}<q^{(n-r-s-2)/2}B$, we conclude
that\[
q^{-n}\sum_{\mathbf{a}\in V}S_{1}(\mathbf{a})S_{2}(\mathbf{a})=q^{-r}N(\Lambda,\sB,q)+O_{n,d}(q^{(n-r-s-2)/2}B).\]

We now estimate $q^{-n}\sum_{\mathbf{a}\in\mathbb{F}_{q}^{n}\setminus V}S_{1}(\mathbf{a})S_{2}(\mathbf{a})$.
Since $\dim Z_{q}\cap\Pi_{\mathbf{a}}<\dim Z_{q}$ for $\mathbf{a}\notin V$,
we obtain from the theorem of Katz (cf. \cite{luo}) that \[
S_{2}(\mathbf{a})\ll_{n,d}q^{(n-r-s+\delta)/2}\]
where $\delta=\dim\Sing(Z_{q}\cap\Pi_{\mathbf{a}})<\dim Z_{q}\in\{-1,0\}$.
As \[
\sum_{\mathbf{a}\in\mathbb{F}_{q}^{n}}|S_{1}(\mathbf{a})|\ll_{n,d}q^{n}(\log q)^{n}\]
(see \cite{luo}), we get that the total contribution to $q^{-n}\sum_{\mathbf{a}\in\mathbb{F}_{q}^{n}\setminus V}S_{1}(\mathbf{a})S_{2}(\mathbf{a})$
from all $\mathbf{a}\in\mathbb{F}_{q}^{n}\setminus V$ where $Z_{q}\cap\Pi_{\mathbf{a}}$
is non-singular is $O_{n,d}(q^{(n-r-s-1)/2}(\log q)^{n})$.

To estimate the contribution from the remaining $\mathbf{a}\in\mathbb{F}_{q}^{n}$,
we use that there exists a form $\Phi\in\mathbb{Z}[y_{1},\ldots,y_{n}]$
of degree $O_{n,d}(1)$ in the dual coordinates $(y_{1},\ldots,y_{n})$
of $(x_{1},\ldots,x_{n})$ such that $\Phi(\mathbf{a})=0$ in $\mathbb{Z}/q\mathbb{Z}$
for all $n$-tuples $\mathbf{a}$ where $Z_{q}\cap\Pi_{\mathbf{a}}$
is singular (cf. Lemma 2 in \cite{marmon}). Hence,\[
\sum_{\substack{\mathbf{a}\in\mathbb{F}_{q}^{n}\\
\Sing(Z_{q}\cap\Pi_{\mathbf{a}})\neq\emptyset}
}|S_{1}(\mathbf{a})|\leq\sum_{\substack{\mathbf{a}\in\mathbb{F}_{q}^{n}\\
\Phi(\mathbf{a})=0}
}|S_{1}(\mathbf{a})|\ll_{n,d}q^{n-1}B(\log q)^{n-1},\]
where the last inequality comes from an argument in \cite{luo}. The
$n$-tuples $\mathbf{a}$ where $Z_{q}\cap\Pi_{\mathbf{a}}$ is singular
will therefore contribute with \[
O_{n,d}(q^{(n-r-s-2)/2}B(\log q)^{n-1})\]
 to $q^{-n}\sum_{\mathbf{a}\in\mathbb{F}_{q}^{n}}S_{1}(\mathbf{a})S_{2}(\mathbf{a})$.
This completes the proof of the lemma.
\end{proof}
For a linear form $L=a_{1}x_{1}+\ldots+a_{n}x_{n}\in\mathbb{Z}[x_{1},\ldots,x_{n}]$,
we will write $\left\Vert L\right\Vert =\sup(|a_{1}|,\ldots,|a_{n}|)$.

\begin{appthm}
\label{thm:appendix}Let $q$ be a prime and $\sB$ be a box in $\mathbb{R}^{n}$
such that each side has length at most $2B<q$. Let $\fs$ be polynomials
in $\mathbb{Z}[x_{1},\ldots,x_{n}],$ $r<n$ with leading forms $F_{1},\ldots,F_{r}$
of degree $\geq2$. Let \begin{gather*}
X=\Spec\mathbb{Z}[x_{1},\ldots,x_{n}]/(f_{1},\ldots,f_{r})\text{ and}\\
Z=\Proj\mathbb{Z}[x_{1},\ldots,x_{n}]/(F_{1},\ldots,F_{r})\end{gather*}
Suppose that $Z_{q}=Z_{\mathbb{F}_{q}}$ is a closed subscheme of
$\mathbb{P}_{\mathbb{F}_{q}}^{n-1}$of codimension $r$ with singular
locus of dimension $s$. Then,\[
N(X,\sB,q)=q^{-r}N(\mathbb{A}_{\mathbb{Z}}^{n},\sB,q)+O_{n,d}(B^{s+1}q^{(n-r-s-2)/2}(B+q^{1/2})(\log q)^{n}),\]
where $d=\max_{i}\deg F_{i}$. 
\end{appthm}
\begin{proof}
It is enough to prove the statement for $q$ greater than some constant
$q_{0}$ depending only on $n$ and $d$, since for $q\ll_{n,d}1$
we have $B\ll_{n,d}1$ and thus, trivially, $N(X,\sB,q)-q^{-r}N(\mathbb{A}_{\mathbb{Z}}^{n},\sB,q)\ll_{n,d}1$.
Thus, assuming that $q$ is large enough, we choose $s+1$ linear
forms $L_{1},\ldots,L_{s+1}\in\mathbb{Z}[x_{1},\ldots,x_{n}]$ such
that $\left\Vert L_{i}\right\Vert =O_{d,n}(1)$ and such that \[
Z_{q}^{i}=\Proj\mathbb{Z}[x_{1},\ldots,x_{n}]/(q,F_{1},\ldots,F_{r},L_{1},\ldots,L_{i})\]
 is a closed subscheme of codimension $r+i$ in $\mathbb{P}_{\mathbb{F}_{q}}^{n-1}$
with singular locus of dimension $s-i$ for $i=1,\ldots,s+1$. Such
forms were used already in \cite{heath-brown} and one gets a proof
of their existence from Lemma 2 in \cite{marmon}.

Let $I=L(\sB\cap\mathbb{Z}^{n})$ for the map $L:\mathbb{Z}^{n}\to\mathbb{Z}^{s+1}$
which sends $\mathbf{b}=(b_{1},\ldots,b_{n})$ to $(L_{1}(\mathbf{b}),\ldots,L_{s+1}(\mathbf{b}))$.
Then $\# I=O_{n,d}(B^{s+1}).$ Moreover, if $\mathbf{c}=(c_{1},\ldots,c_{s+1})\in\mathbb{Z}^{s+1}$,
then we may apply Lemma \ref{lem:appendix} to the affine subscheme
$X_{\mathbf{c}}$ of $\mathbb{A}_{\mathbb{Z}}^{n}$ defined by $(f_{1},\ldots,f_{r},L_{1}-c_{1},\ldots,L_{s+1}-c_{s+1})$
and conclude that\[
N(X_{\mathbf{c}},\sB,q)=q^{-r}N(\Lambda_{\mathbf{c}},\sB,q)+O_{n,d}(q^{(n-r-s-2)/2}(B+q^{1/2})(\log q)^{n})\]
for $\Lambda_{\mathbf{c}}=\Spec\mathbb{Z}[x_{1},\ldots,x_{n}]/(L_{1}-c_{1},\ldots,L_{s+1}-c_{s+1})$.
If we sum over all $\mathbf{c}=(c_{1},\ldots,c_{s+1})\in I$, then
we get the desired asymptotic formula for $N(X,\sB,q)$. This finishes
the proof.
\end{proof}
\begin{rem*}
Note that $q^{-r}N(\mathbb{A}_{\mathbb{Z}}^{n},\sB,q)=q^{-r}\#(\sB\cap\mathbb{Z}^{n})$,
since different elements in $\sB\cap\mathbb{Z}^{n}$ are non-congruent
(mod $q$) by the assumption on $\sB$.
\end{rem*}

\end{document}